\documentclass[a4paper,11pt]{article}

\usepackage{latexsym}
\newtheorem{thm}{Theorem}[section]

\newtheorem{cor}[thm]{Corollary}
\newtheorem{prop}[thm]{Proposition}

\newtheorem{rem}[thm]{Remark}
\input epsf

\title{Counting Triangulations of Configurations}
\author{Roland Bacher
}
\begin{document}
\maketitle
%\par fichier etoile1.tex dans recherche/triangulations

%Identities (Orthogonal polys for $p_n$),

\section{Introduction}                              

Given a finite subset ${\mathcal C}\subset {\mathbf R}^2$,
calculating the number of triangulations of its convex hull 
$\hbox{Conv}({\mathcal C})$  using only Euclidean triangles
with vertices in ${\mathcal C}$ seems to be difficult and has attracted
some interest, both from an algorithmic and a theoretical point of view,
see for instance \cite{A1}, \cite{A2}, \cite{AHN}, \cite{AK}, \cite{AAK},
\cite{An}, \cite{KZ}, \cite{S}, \cite{SS}.

The aim of this paper is to describe a class of configurations,
{\it convex near-gons}, for which this problem can sometimes be solved in a
satisfactory way. Loosely speaking, a convex near-gon is an infinitesimal
perturbation of a {\it weighted convex polygon}, a convex 
polygon with edges subdivided
by additional points according to weights. 
Our main result shows that the
triangulation polynomial enumerating all triangulations of a convex 
near-polygon
is defined in a straightforward way in terms of
edge-polynomials associated to the \lq\lq perturbed'' edges 
of such a convex near-gone. These polynomials 
are difficult to compute in general except in a few special cases.
We present a few algorithms related to them.
One of these algorithms is a slightly more sophisticated version of
an algorithm by Kaibel and Ziegler described in \cite{KZ}
and yields also a general purpose 
algorithm (unfortunately of exponential complexity), for computing
arbitrary triangulation polynomials. 
This algorithm, based on a transfer matrix, is fairly simple and it would
be interesting to compare its performance with existing 
algorithms, like for instance the algorithm of Aichholzer described
in \cite{A1}.

\begin{figure}[h]\label{exple1}
\epsfysize=4cm
\epsfbox{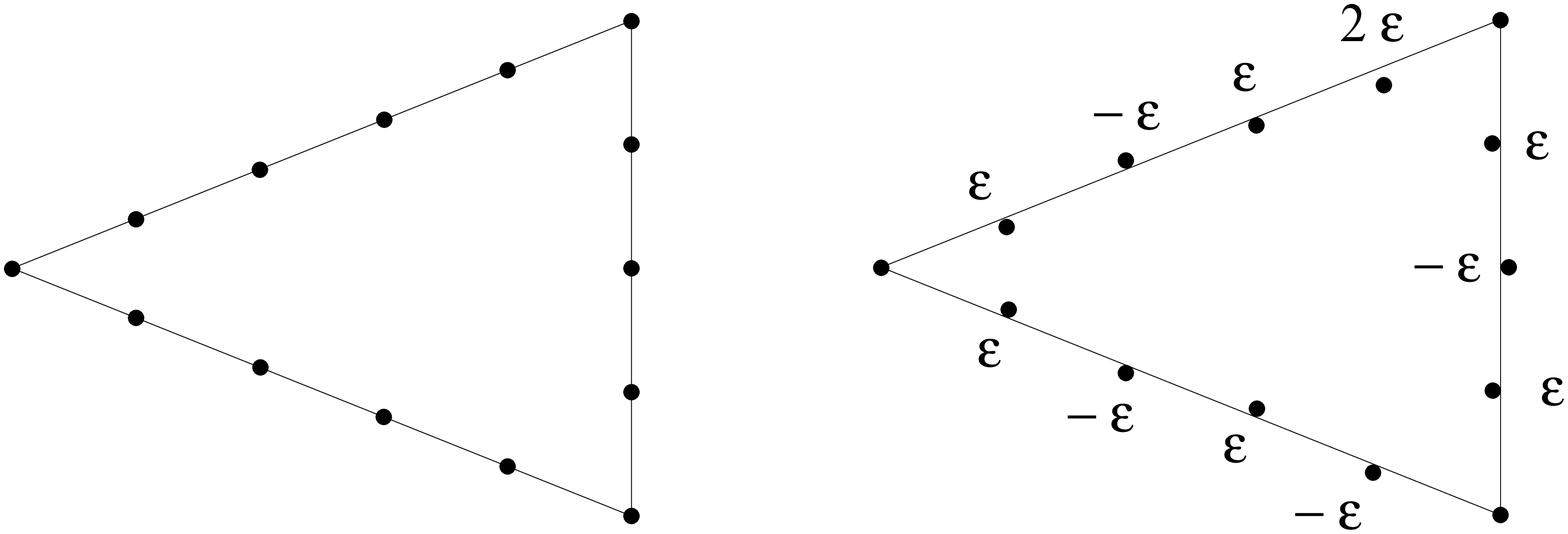}
\caption{A weighted triangle and an associated generic perturbation}
\end{figure}  
 
\begin{figure}[h]\label{expleconfig}
\epsfysize=7cm
\epsfbox{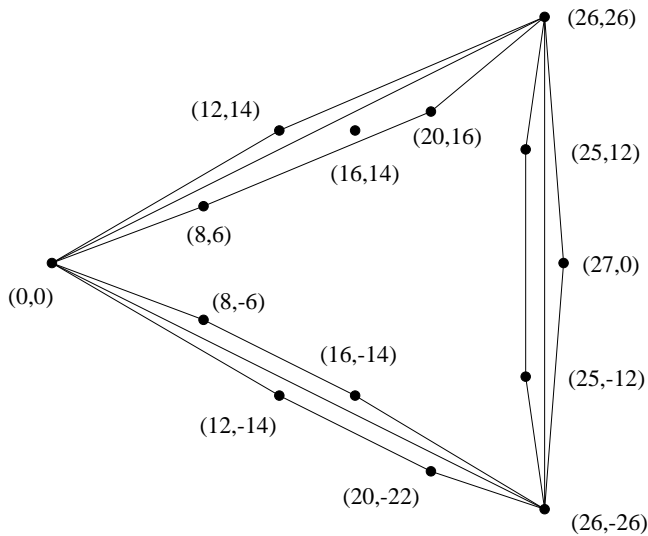}
\caption{An integral representant of a convex near-polygon}
\end{figure}  

Figure \ref{exple1} represents a weighted convex triangle and an 
infinitesimal generic perturbation of it. Figure 2
shows an integral realisation isotopic to this infinitesimal 
perturbation.

The complete triangulation polynomial of the weighted triangle $\Delta$
depicted on the left side of Figure \ref{exple1} equals
$$\begin{array}{c}
\displaystyle 
901\,{s}^{14}+4825\,{s}^{13}+11734\,{s}^{12}+17130\,{s}^{11}+16710\,{s
}^{10}+11466\,{s}^{9}
\\
\displaystyle
+5670\,{s}^{8}+2034\,{s}^{7}+525\,{s}^{6}+95\,{s}
^{5}+11\,{s}^{4}+{s}^{3}
\quad
\end{array}$$
showing that there are for instance $17130$ triangulations
of $\Delta$ using $11$ vertices among the $14$ vertices of $\Delta$.

This paper contains all details (except for some straightforward
calculations involving polynomials where computers exceed human
possibilities by orders of magnitudes) needed 
for computing the complete triangulation polynomial
$$\begin{array}{l}
\displaystyle 
194939\,{s}^{14}+338669\,{s}^{13}+263615\,{s}^{12}+119944\,{s}^{11}\\
\displaystyle 
\quad +34773\,{s}^{10}+6522\,{s}^{9}+748\,{s}^{8}+42\,{s}^{7}\end{array}$$
of the near-triangle represented in Figure 2.

In the sequel of this paper, we recall first a few generalities
concerning finite configurations of ${\mathbf R}^2$.

We define then weighted convex polygons and give a formula for their
complete triangulation polynomial in terms of their
edge-weights.

After introducing convex near-gons, near-edges and edge-polynomials
we state the main result of this paper expressing the complete
triangulation polynomial of a convex near-polygon in terms of the
edge-polynomials of its edges.

We list then a few usefull facts and give some data
concerning edge-polynomials. 
 
Finally, we describe a transfer-matrix approach for computing 
edge-polynomials (and/or solving enumerative problems 
concerning triangulations).

We close this introduction by describing the following notation,
suggested by the fact that polynomials and formal
power series are each other's algebraic dual.

\noindent{\bf Notations.} Given a polynomial $p(x)=\sum_{i=0}^d \alpha_i\ x^i
\in {\mathbf k}[x]$
and a formal power series $s(x)=\sum_{i=0}^\infty \beta_i\ x^i
\in {\mathbf k}[[x]]$ (over any commutative field or 
ring ${\mathbf k}$) we denote by
$$\langle p,s\rangle_x=\langle p(x),s(x)\rangle_x=\sum_{i=0}^d
\alpha_i\ \beta_i$$
the obvious linear pairing. 

This notation makes sense for polynomials and formal power 
series over several variables and we have
$$\langle p,s\ t\rangle_{x,y}=\langle\ \langle p,s\rangle_x\ ,t\rangle_y
=\langle\ \langle p,t\rangle_y\ ,s\rangle_x$$
for a polynomial $p\in{\mathbf k}[x,y]$ and formal power series
$s\in{\mathbf k}[[x]],\ t\in{\mathbf k}[[y]]$.

Throughout the rest of this paper, $p_n$ denotes always the polynomial 
$$p_n=\sum_{k=0}^{\lfloor n/2\rfloor}(-1)^k{n-k\choose k}t^{n-k}\in {
\mathbf Z}[t]$$
and $C_n={2n\choose n}/(n+1)$ stands for the $n-$th Catalan number.

%%%%%%%%%%%%%%%%%%%%%%%%%%%%%%%%%%%%%%%%%%%%%%%%%%%%%%%%%%%%%%%%%%%%

\section{Configurations}

A {\it configuration} of the oriented plane
${\mathbf R}^2$ is a finite set ${\mathcal C}=
\{P_1,\dots,P_n\}\subset {\mathbf R}^2$ of $n$ distinct points.
Two configurations ${\mathcal C}=\{P_1\dots,P_n\}$ and
${\mathcal C}'=\{P_1',\dots,P_n'\}$ of $n$ points are {\it (combinatorially)
equivalent} or {\it isomorphic} if there exists a bijection
$\varphi:{\mathcal C}\longrightarrow {\mathcal C}'$ such that
$$\hbox{sign}(\det(P_j-P_i,P_k-P_i))=
\hbox{sign}(\det(\varphi(P_j)-\varphi(P_i),\varphi(P_k)-\varphi(P_i)))$$
for all $1\leq i<j<k\leq n$
where
$$\hbox{sign}(x)=\left\{\begin{array}{rcl}
1&&\hbox{if }x>0\\
0&&\hbox{if }x=0\\
-1&\quad&\hbox{if }x<0\end{array}\right.\ .$$
Two configurations ${\mathcal C}=\{P_1,\dots,P_n\}$ and
${\mathcal C}'=\{P_1',\dots,P_n'\}$ of $n$ points are {\it isotopic} 
if there exists a continuous path $t\longmapsto {\mathcal C}(t)\in
\left({\mathbf R}^2\right)^n$
of equivalent configurations such that ${\mathcal C}(-1)={\mathcal C}$
and
${\mathcal C}(1)={\mathcal C}'$. Isotopic configurations
are equivalent.

A configuration ${\mathcal C}\subset{\mathbf R}^2$ is {\it generic} if
three distinct points of ${\mathcal C}$ are never collinear.

Given a configuration ${\mathcal C}$, we denote by $\hbox{Conv}({\mathcal
C})\subset {\mathbf R}^2$ its convex hull. An element
$P\in {\mathcal C}$ is {\it extremal} if $\hbox{Conv}({\mathcal C}
\setminus\{P\})\not= \hbox{Conv}({\mathcal C})$.
We denote by $\hbox{Extr}({\mathcal C})
\subset {\mathcal C}$ the subset of extremal vertices. They
are the vertices of the convex polygone $\hbox{Conv}({\mathcal C})$.

A {\it triangulation} of a finite configuration ${\mathcal C}$ is a 
finite set ${\mathcal T}=\{\Delta_1,\dots,\Delta_k\}$
of Euclidean triangles with vertices in ${\mathcal C}$
such that $\hbox{Conv}({\mathcal C})=\cup_{\Delta_i\in{\mathcal T}}
\Delta_i$ and $\Delta_i\cap \Delta_j$ is either empty, or a common
vertex, or a common edge of two distinct triangles $\Delta_i,\Delta_j\in
{\mathcal T}$.

The {\it complete triangulation polynomial} ${\overline p}({\mathcal C})$ 
is the polynomial ${\overline p}({\mathcal C})=\sum \tau_k({\mathcal
C})\ s^k\in{\mathbf Z}[s]$ 
whose coefficient $\tau_k({\mathcal C})$ equals the number
of triangulations of $\hbox{Conv}({\mathcal C})$ using exactly $k$ vertices of 
${\mathcal C}$. The number of {\it maximal triangulations}
using all vertices of ${\mathcal C}$ is given by the leading coefficient of 
${\overline p}({\mathcal C})$. The lowest non-zero coefficient
of  ${\overline p}({\mathcal C})$ is the $(m-2)-$th Catalan number
$C_{m-2}={2(m-2)\choose m-2}/(m-1)$ where $m=\sharp(\hbox{Extr}
({\mathcal C}))$ is the number of vertices of the polygone
$\hbox{Conv}({\mathcal C})$, see for instance Exercice 6.19 of
\cite{St}. It is easy to check that the complete
triangulation polynomial ${\overline p}({\mathcal C})$ depends only
of the equivalence class of a finite configuration ${\mathcal C}$.

\begin{rem} One can also consider the triangulation
polynomial defined by
$$\sum_{k_0,k_1}\tau_{k_0,k_1}({\mathcal C})s_0^{k_0}s_1^{k_1}$$
counting the number of triangulations using $k_0$ vertices
and $k_1$ edges. The number $k_2$ of triangles can then be recovered 
using the Euler characteristic $k_0-k_1+k_2=1$ of a compact, simply 
connected triangulated polygonal region in ${\mathbf R}^2$.
This more general polynomial yields the same information as the
complete polynomial considered above except if the boundary
$\partial(\hbox{Conv}({\mathcal C}))$ contains points of 
${\mathcal C}$ which are not extremal. Most of the results
and algorithms 
of this paper can easily be modified in order to deal with this
more general polynomial. For clarity and concision we stick to the
simpler version $\overline p({\mathcal C})$ defined above.
\end{rem}
%%%%%%%%%%%%%%%%%%%%%%%%%%%%%%%%%%%%%%%%%%%%%%%%%%%%%%%%%%%%%%%

\section{Weighted convex polygons}

This section describes weighted convex polygons,
a particular set of finite planar configurations,
and gives formulae for their number of maximal
triangulations and for their complete triangulation polynomials.

\begin{figure}[h]\label{weighted-pentagon15234}
\epsfysize=3cm
\epsfbox{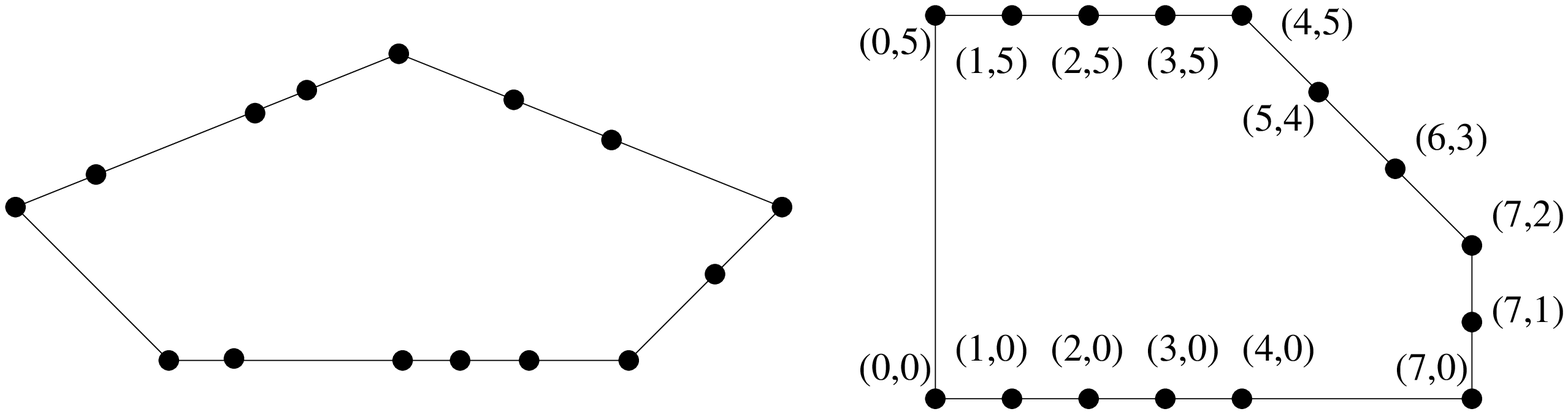}
\caption{
Two weighted pentagons with edge-weights $1,5,2,3,4$.}
\end{figure}

For $l\geq 3$ and natural numbers
$a_1,\dots,a_l\geq 1$ we denote by $P(a_1\cdot a_2\cdots a_l)
\subset{\mathbf R}^2$
a strictly convex polygon with $l$ (counterclockwise) 
cyclically ordered edges
subdivided respectively by $(a_1-1),(a_2-1),\dots,(a_l-1)$
additional points. We call such a subdivided polygon a
{\it weighted polygon} with (cyclical) weights $a_1,\dots,a_l$. 
Figure \ref{weighted-pentagon15234} displays two realisations of
a weighted pentagon $P(1\cdot 5\cdot 
2\cdot 3\cdot 4)$.

The configurations given by the $(a_1+\dots +a_l)=(b_1+\dots +b_l)$ 
vertices of two weighted convex polygons $P(a_1\cdots a_l)$ and
$P(b_1\cdots b_l)$ having the same number of edges and the same
weights up to cyclical permutation, are isotopic and
their sets of triangulations have the same combinatorics. 
We denote by $\tau_k(a_1\cdots a_l)$ the number
of triangulations which use $k$ vertices of such a
weighted convex polygon $P(a_1\cdots a_l)$. 
Such triangulations exist of course only if
$l\leq k\leq \sum_{i=1}^l a_i$.
 
For $n\geq 1$ we 
introduce the {\it maximal edge-polynomial} 
$p_n\in{\mathbf Z}[t]$ by setting
$$p_n=\sum_{k=0}^{\lfloor n/2\rfloor}(-1)^k{n-k\choose k}t^{n-k}\in {
\mathbf Z}[t]$$
and the {\it complete edge-polynomial} ${\overline p}_n\in
\left({\mathbf Z}[s]\right)[t]$ with indeterminate $t$ over the
ring ${\mathbf Z}[s]$ which is defined as
$${\overline p}_n=\sum_{k=1}^n {n-1\choose k-1}p_k(t)s^k\in
\left({\mathbf Z}[s]\right)[t]\ .$$

The first few maximal edge-polynomials are
$$\begin{array}{ll}
\displaystyle p_1=t&
\displaystyle p_5=t^5-4t^4+3t^3\\
\displaystyle p_2=t^2-t&
\displaystyle p_6=t^6-5t^5+6t^4-t^3\\
\displaystyle p_3=t^3-2t^2&
\displaystyle p_7=t^7-6t^6+10t^5-4t^3\\
\displaystyle p_4=t^4-3t^3+t^2\qquad&
\displaystyle p_8=t^8-7t^7+15t^6-10t^5+t^4\\
\end{array}$$
and the first few complete edge-polynomials are
$$\begin{array}{l}
\displaystyle {\overline p}_1=p_1\ s=t\ s\ ,\\
\displaystyle {\overline p}_2=p_2\ s^2+p_1\ s=(t^2-t)\ s^2+t\ s\ ,\\
\displaystyle {\overline p}_3=p_3\ s^3+2p_2\ s^2+p_1\ s=(t^3-2t^2)\ s^3+
2(t^2-t)\ s^2+t\ s\ ,\\
\displaystyle {\overline p}_4=p_4\ s^4+3p_3\ s^3+3p_2\ s^2+p_1\ s\\
\displaystyle \qquad =(t^4-3t^3+t^2)\ s^4+3(t^3-2t^2)\ s^3+
3(t^2-t)\ s^2+t\ s\ .\end{array}$$

\begin{thm} \label{maxformulaweighted} 
The number $\tau_{\hbox{max}}(a_1\cdot a_2\cdots a_l)$ 
of maximal triangulations of the weighted convex polygon
$P(a_1\cdot a_2\cdots a_l)$ is given by
$$\langle \prod_{j=1}^l {p}_{a_j}(t),\sum_{n=2}^\infty C_{n-2}\ t^n
\rangle_t$$
where $C_n={2n\choose n}/(n+1)$  is the $n-$th Catalan number.
\end{thm}

\begin{cor} \label{completeformulaweighted} The complete triangulation polynomial 
$\sum_k\tau_k(a_1\cdots a_l)s^k$ of a weighted convex
polygon $P(a_1\cdot a_2\cdots a_l)$ is given by
$$\langle \prod_{j=1}^l {\overline p}_{a_j}(t),\sum_{n=2}^\infty
C_{n-2}\ t^n\rangle_t\ .$$
\end{cor}

\begin{cor} \label{orderindepfor weighted}
The complete triangulation polynomial (and hence the number 
of maximal triangulations) of
a convex weighted polygone $P(a_1\cdots a_l)$ depends only on the multi-set
(set having perhaps multiple elements)
$\{a_1,\dots,a_l\}$ defined by the weights 
and not on their particular cyclic order.
\end{cor}

\noindent{\bf Example.} The two weighted convex polygons $P(1,5,2,3,4)$
of Figure \ref{weighted-pentagon15234} have 
$$\begin{array}{cl}\displaystyle 
& \langle p_1\ p_5\ p_2\ p_3\ p_4,\sum_{n=2}^\infty  C_{n-2}\ t^n
\rangle_t\\ \displaystyle 
=&\langle t(t^5-4t^4+3t^3)(t^2-t)(t^3-2t^2)(t^4-3t^3+t^2),\sum_{n=2}^\infty
C_{n-2}\ t^n\rangle_t\\ \displaystyle
=&\langle t^{15}-10t^{14}+39t^{13}-75t^{12}
+74t^{11}-35t^{10}+6t^9,\sum_{n=2}^\infty C_{n-2}\ 
t^n\rangle_t\\ \displaystyle
=&C_{13}-10C_{12}+39C_{11}-75C_{10}+74C_{9}-35C_8+6C_7\\ \displaystyle
=&7429000-10\cdot 208012+39\cdot 58786-75\cdot 16796\\
\displaystyle& \quad +74\cdot 4863-35\cdot
1430+6\cdot 429\\ \displaystyle
=&8046\end{array}$$
maximal triangulations. 

Their complete triangulation polynomial 
$\sum_k \tau_k(P(1\cdot 5\cdot 2\cdot 3\cdot 4))s^k$ equals
$$\begin{array}{cl}\displaystyle 
& \langle {\overline p}_1\ {\overline p}_5\ {\overline p}_2\ 
{\overline p}_3\ {\overline p}_4,\sum_{n=2}^\infty C_{n-2}t^n
\rangle_t\\ \displaystyle =&
8046\,{s}^{15}+37250\,{s}^{14}+77467\,{s}^{13}+95364\,{s}^{12}+77048\,
{s}^{11}\\
\displaystyle &\quad +42776\,{s}^{10}+16584\,{s}^{9}+4460\,{s}^{8}+805\,{s}^{7}+90
\,{s}^{6}+5\,{s}^{5}\ .
\end{array}$$

\begin{figure}[h]\label{digon}
\epsfysize=1.8cm
\epsfbox{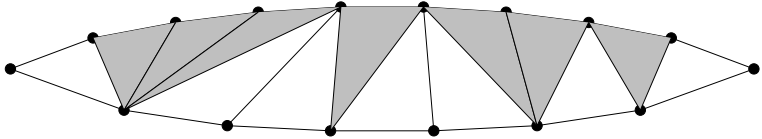}
\caption{A maximal triangulation with 7 (shaded) up-triangles of 
the digon $P(9\cdot 7)$.}
\end{figure}

\begin{rem} \label{digons}
Theorem \ref{maxformulaweighted} and Corollary 
\ref{completeformulaweighted} remain valid for weighted \lq\lq di\-gons'',
weighted \lq\lq convex'' polygons $P(a_1\cdot a_2)$ having only two 
edges. Their \lq\lq sides'' or \lq\lq edges'' are in this case not straight
lines but piecewise linear paths bending sligthly outwards in order to 
enclose a non-empty interior if $a_1\ a_2\geq 2$ 
(see Figure 4 for an example). 
The weighted digon $P(1,1)$ has by convention one triangulation.
$P(0,a)=P(a,0)$ has 1 triangulation for $a=0$ and none otherwise.
In all the remaining cases, only \lq\lq 
triangulations'' with all triangles meeting both edges
are allowed. The number of such triangulations of $P(a_1\cdot a_2)$
equals ${a_1+a_2-4\choose a_1-2}$: Indeed, every triangle except the
first and the last one is either up (with two vertices on the upper edge)
or down (with two vertices on the lower edge) and there are exactly $a_1-2$
such triangles which are up, cf. Figure 4 for an example.
\end{rem}

\subsection{Proofs for weighted convex polygons}

Let $P(a_1^{\alpha_1}\cdot a_2^{\alpha_2}\cdots a_l^{\alpha_l})$ denote
the weighted convex polygon having $\alpha_1$ edges of weight $a_1$ 
followed by $\alpha_2$ edges of weight $a_2$ etc. We denote by
$\tau_{\hbox{max}}(a_1^{\alpha_1}\cdot a_2^{\alpha_2}\cdots a_l^{\alpha_l})$
the number of maximal triangulations of 
$P(a_1^{\alpha_1}\cdots a_l^{\alpha_l})$. 

The main ingredient for proving
Theorem \ref{maxformulaweighted} is the following illustration 
of the Inclusion-Exclusion Principle (cf. Chapter 2 of of \cite{St1}).

\begin{prop} \label{inclusion-exclusion}
(Inclusion-Exclusion Principle.) We have
$$\tau_{\hbox{max}}(a_1\cdot a_2\cdots a_l)=\sum_{k=0}^{\lfloor a_1/2\rfloor}
(-1)^k{a_1-k\choose k}\tau_{\hbox{max}}(1^{a_1-k}\cdot a_2\cdots a_l)\ .$$
\end{prop}

\noindent{\bf Proof.} Enumerate
the sides of $P(1^{a_1}\cdot a_2\cdots a_l)$ cyclically such that the 
first $a_1$ sides are of weight $1$ and correspond to the initial
factor $1^{a_1}$ of the product $1^{a_1}\cdot a_2\cdots a_l$.
Let $T$ be a maximal triangulation of the
weighted convex polygon $P(1^{a_1}\cdot a_2 \cdots a_l)$. An {\it 
initial boundary triangle}
of $T$ is a triangle of $T$ having
two edges contained among the first $a_1$ sides (having weight $1$) 
of the convex weighted polygon $P(1^{a_1}\cdot a_2 \cdots a_l)$.  
A subset $\{\delta_1,\dots,\delta_k\},\ 0\leq k\leq \lfloor a_1/2\rfloor$
of $k$ initial boundary triangles in $T$ is called a {\it
$\Delta_k-$decoration} of
$T$. Starving to death the $k$ triangles $\delta_1,\dots,\delta_k$
of the $\Delta_k-$decorated maximal
triangulation $(T,\{\delta_1,\dots,\delta_k\})$,
we get a maximal triangulation $\tilde T$ of $P(1^{a_1-k}\cdot a_2\cdots a_l)$
together with a graveyard $\{e_1,\dots,e_k\}$ of $k$ marked edges in memory 
of the deceased triangles. We call the pair $(\tilde T,\{e_1,\dots,e_k\})$
a maximal {\it $E_k-$decorated} triangulation 
of $P(1^{a_1-k}\cdot a_2\cdots a_l)$. Figure 5
%\ref{2decorated} 
displays a $\Delta_2-$decorated  maximal triangulation of 
$P(1^7\cdot 3\cdot 2)$ and the corresponding 
$E_2-$decorated maximal triangulation of $P(1^5\cdot 3\cdot 2)$.

\begin{figure}[h]\label{2decorated}
\epsfysize=3.2cm
\epsfbox{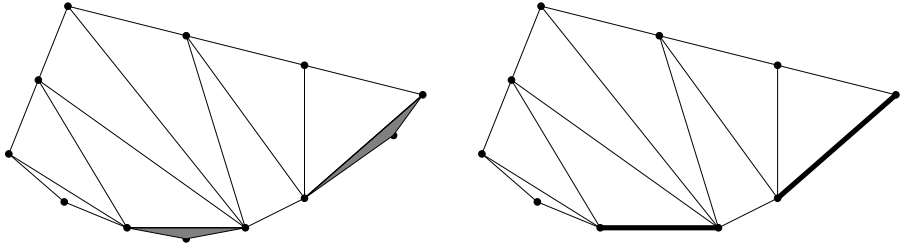}
\caption{Decorated  maximal triangulations of $P(1^7\cdot 3\cdot 2)$ and 
$P(1^5\cdot 3\cdot 2)$}
\end{figure}  

Obviously, $\Delta_k-$decorated maximal triangulations of 
$P(1^{a_1}\cdot a_2\cdots a_l)$ and $E_k-$decorated maximal triangulations
of $P(1^{a_1-k}\cdot a_2\cdots a_l)$ are in bijection and since any 
maximal triangulation of $P(1^{a_1-k}\cdot a_2\cdots a_l)$ has 
${a_1-k\choose k}$ different $E_k-$decorations, there exist exactly
$${a_1-k\choose k}\tau_{\hbox{max}}(1^{a_1-k}\cdot a_2\cdots a_l)$$
maximal $\Delta_k-$decorated triangulations of 
$P(1^{a_1}\cdot a_2\cdots a_l)$.

\begin{figure}[h]\label{straightening}
\epsfysize=1.8cm
\epsfbox{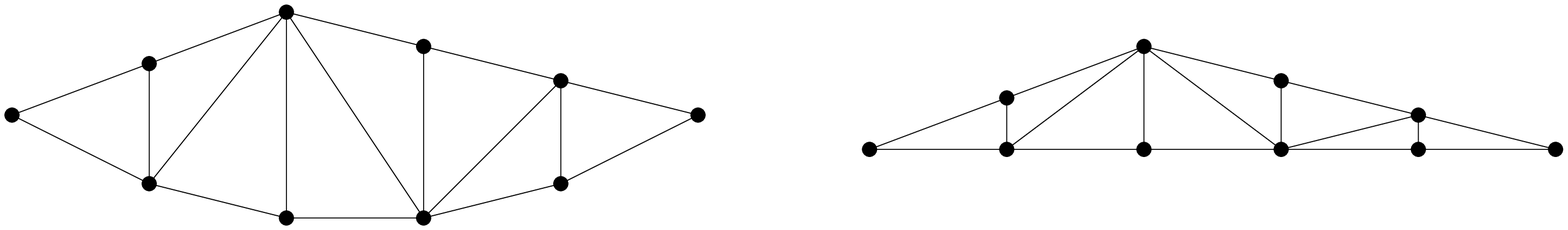}
\caption{Straightening a maximal triangulation without initial 
boundary triangles}
\end{figure}

Consider now a maximal triangulation $T$ of $P(1^{a_1}\cdot a_2\cdots a_l)$
having exactly $s$ initial boundary triangles. The triangulation $T$ 
gives rise to ${s\choose k}$ different $\Delta_k-$decorated maximal 
triangulations of $P(1^{a_1}\cdot a_2\cdots a_l)$ and yields
a contribution of 
$$\sum_{k=0}^s{s\choose k}(-1)^k=\left\lbrace\begin{array}{ll}
\displaystyle 1\quad &\hbox{if }s=0\\
\displaystyle 0&\hbox{otherwise}\end{array}\right.$$
to the alternating sum 
$$\sum_{k=0}^{\lfloor a_1/2\rfloor}
(-1)^k{a_1-k\choose k}\tau_{\hbox{max}}(1^{a_1-k}\ a_2\cdots a_l)$$
counting $\Delta_k-$decorated maximal triangulations of 
$P(1^{a_1}\cdot a_2\cdots a_l)$ with sign $(-1)^k$. 
This alternating sum counts thus exactly the number of 
maximal triangulations 
of $P(1^{a_1}\cdot a_2\cdots a_l)$ without initial
boundary triangles. Such triangulations are in bijection
with triangulations of $P(a_1\cdot a_2\cdots a_l)$:
straight out the first $a_1$ edges of weight $1$ in 
$P(1^{a_1}\cdot a_2\cdots a_l)$ into a unique first edge of weight
$a_1$ of $P(a_1\cdots a_l)$. Figure \ref{straightening} 
illustrates this: Its left side displays a maximal triangulation
without initial boundary triangles of $P(1^5\cdot 3\cdot 2 )$
and its right side shows the corresponding straightened
maximal triangulation of $P(5\cdot 3\cdot 2)$. \hfill $\Box$

\noindent{\bf Proof of Theorem \ref{maxformulaweighted}}
Given a polynomial $q(t)=\sum_{i=0}^d \alpha_i\ t^i$ we set
$$\tau(q(t);\ a_2\cdots a_l)=\sum_{i=0}^d \alpha_i\  \tau_{\hbox{max}}
(1^i\cdot a_2\cdots a_l)\ .$$
Since the polynomial
$$p_{a_i}(t)=\sum_{k=0}^{\lfloor a_i/2\rfloor} (-1)^k{n-k\choose k}t^{n-k}$$
mimicks the formula of Proposition \ref{inclusion-exclusion},
we get by iterated application of Proposition \ref{inclusion-exclusion}
$$\tau_{\hbox{max}}(a_1\cdot a_2\cdots a_l)=
\tau(\prod_{i=1}^l p_{a_i}(t);\emptyset)=
\sum_j\gamma_j\  \tau_{\hbox{max}}(1^j)$$
where $\sum_j \gamma_j t^j=\prod_{i=1}^l p_{a_i}(t)$.

Since $P(1^n)$ is a convex polygon with $n$ edges (of weight
$1$) and $n$ extremal vertices, the number $\tau_{\hbox{max}}(1^n)$ of
maximal triangulations of $P(1^n)$ is given by the $(n-2)-$th Catalan number
$C_{n-2}={2(n-2)\choose n-2}/(n-1)$, cf. for instance 
Exercice 6.19 of \cite{St}. 
Hence the result.
\hfill$\Box$

\noindent{\bf Proof of Corollary \ref{completeformulaweighted}.} Non-maximal
triangulations of $P(a_1\cdots a_l)$ are in bijection with maximal
triangulations of $P(b_1\cdots b_l)$ where $1\leq b_i\leq a_i$ 
and where the weighted polygons $P(b_1\cdots b_l)$ are realised by
choosing in all ${a_i-1\choose b_i-1}$ possible ways
subsets of $(b_i-1)$ points among the $(a_i-1)$ interior points of the 
$i-$th edge having weight $a_i$ of $P(a_1\cdots a_l)$.

The complete triangulation polynomial of $P(a_1\cdots a_l)$ is thus
given by
$$\begin{array}{l}
\displaystyle \sum_{1\leq b_i\leq a_i}\left(\prod_{j=1}^l {a_i-1\choose b_i-1}
\right)\big\langle\prod_{j=1}^lp_{b_j}(t),\sum_{n=2}^\infty C_{n-2}\ t^n
\big\rangle_t\quad s^{\sum_{j=1}^l b_j}\\ \displaystyle 
\qquad =\big\langle\prod_{j=1}^l\sum_{b_j=1}^{a_j}{a_j-1\choose b_j-1}
p_{b_j}(t)s^{b_j},\sum_{n=2}^\infty C_{n-2}\ t^n\big\rangle_t
\\ \displaystyle 
\qquad =\big\langle\prod_{j=1}^l{\overline p}_{a_j},
\sum_{n=2}^\infty C_{n-2}\ t^n\big\rangle_t
\end{array}
$$
which proves the result. \hfill $\Box$

Corollary \ref{orderindepfor weighted} is obvious.

%%%%%%%%%%%%%%%%%%%%%%%%%%%%%%%%%%%%%%%%%%%%%%%%%%%%%%%%%%%%%%%%%%%%%%%%

\section{Near-gons}

An $n-${\it near-edge} $E$ is a sequence of $(n+1)-$points
$$(P_0,\dots,P_n)=(\left(\begin{array}{c}x_0\\y_0\end{array}\right),
\left(\begin{array}{c}x_1\\y_1\end{array}\right),\dots,
\left(\begin{array}{c}x_{n-1}\\y_{n-1}\end{array}\right),
\left(\begin{array}{c}x_{n}\\y_n\end{array}\right))$$
in ${\mathbf R}^2$ such that $x_0<x_1<\dots<x_{n-1}<x_{n}$. We 
consider the complete order $P_0<P_1<P_2<\dots,P_n$ on $E$ and 
call $P_0$, respectively $P_n$, the {\it initial}, respectively {\it final},
vertex of the near-edge $E$. We denote a
near-edge $E$ either by the sequence $E=(P_0,\dots,P_n)$
of its points or by the  $2\times (n+1)$ matrix
$$E=\left(\begin{array}{ccc}x_0&\dots&x_n\\ y_0&\dots&y_n\end{array}\right)$$
satisfying $x_0<x_1\dots<x_n$ whose columns contain 
the coordinates
$\left(\begin{array}{c}x_i\\y_i\end{array}\right)$ of $P_i$. We endow $E$
with the total order $P_0<P_1<\dots<P_n$.

Two $n-$near-edges $E=\{P_0,\dots,P_{n}\}$ and $E'=\{P'_0,\dots,P'_{n}\}$
are {\it equivalent} if 
$$\hbox{sign}\left(\det(P_j-P_i,P_k-P_i)\right)=
\hbox{sign}\left(\det(P'_j-P'_i,P'_k-P'_i)\right)$$
where $0\leq i<j<k\leq n$ and where $\hbox{sign}(x)\in\{\pm 1,0\}$
is defined by 
$$\hbox{sign}(x)=\left\{\begin{array}{rcl}
1&&\hbox{if }x>0\\
0&&\hbox{if }x=0\\
-1&\quad&\hbox{if }x<0\end{array}\right.\ .$$

A continuous path $t\longmapsto E(t),\ t\in[-1,1]$ 
of equivalent near-edges is an
{\it isotopy}. Two near-edges $E(-1)$ and $E(1)$ 
joined by an isotopy are {\it isotopic}.
Isotopic near-edges are of course always equivalent.

A near-edge is {\it generic} it the underlying set of points
is a generic configuration of ${\mathbf R}^2$, i.e. is without three 
collinear points.

Given a near-edge $E=\big(P_i=\left(\begin{array}{c}x_i\\y_i\end{array}\right)
\big)_{i=0,\dots,n}$ and $\epsilon\in{\mathbf R}$, 
we denote by $E^\epsilon$ the near-edge
$$E^\epsilon=\big(P_i(\epsilon)=\left(\begin{array}{c}x_i\\\epsilon
\ y_i\end{array}\right)=\left(\begin{array}{cc}1&0\\0&\epsilon\end{array}
\right)P_i\big)_{i=0,\dots,n}\ .$$
For $\epsilon>0$, the near-edges $E$ and $E^\epsilon$ are isotopic and thus
equivalent.

Let $P$ be a convex polygon with $l$ extremal vertices
$\{V_0,V_1,V_2,\dots,V_{l-1},V_l=V_0\}$
appearing in counterclockwise order around its boundary $\partial P$.

Given a sequence $E_1^{\epsilon_1},\dots,E_l^{\epsilon_l}$
where $E_i=(P_{i,0},P_{i,1},\dots,P_{i,n_i})$ 
is a $n_i-$near-edges and $\epsilon_i\in{\mathbf R}$ 
a real number, we denote by
$$P(E_1^{\epsilon_1}\cdot E_2^{\epsilon_2}\cdots E_l^{\epsilon_l})$$
the unique configuration obtained by gluing the 
$n_i-$near-edge $E_i^{\epsilon_i}$, rescaled by a suitable 
orientation-preserving similitude, along the oriented
edge of $P$ which starts at $V_{i-1}$ and ends at $V_i$. More precisely,
the gluing similitude $\varphi_i$ is the unique orientation-preserving 
similitude of ${\mathbf R}^2$ such that 
$$\varphi_i(P_{i,0}(\epsilon_i))=V_{i-1}\hbox{ and }\varphi_i(P_{i,n_i}
(\epsilon_i))=V_{i}\ .$$
The configuration $P(E_1^{\epsilon_1}\cdot E_2^{\epsilon_2}\cdots 
E_l^{\epsilon_l})$ is now the set of points
$\cup_{i=1}^l \varphi_i(E_i^{\epsilon_i}) \subset {\mathbf R}^2$.

We have the following result which we state
without proof.

\begin{prop} \label{near_well-defined}
(i) The configurations
$P(E_1^{\epsilon_1}\cdots E_l^{\epsilon_l})$ are isotopic for 
$0<\epsilon_i$ small enough.

\ \ (ii) Given  a second convex polytope $P'$ having $l$
vertices $V'_0,V'_1,\dots,V'_l=V'_0$ in counterclockwise order, 
the configurations
$$P(E_1^{\epsilon_1}\cdots E_l^{\epsilon_l})\hbox{ and }
P'(E_1^{\epsilon_1}\cdots E_l^{\epsilon_l})$$
associated to $P$ and $P'$ are isotopic for $0< \epsilon_i$ 
small enough.

\ \ (iii) Given pairs $(E_i,E'_i)$ of equivalent near-edges,
the configurations
$$P(E_1^{\epsilon_1}\cdots E_l^{\epsilon_l})\hbox{ and }
P({E'_1}^{\epsilon_1}\cdots {E'_l}^{\epsilon_l})$$
are equivalent for $0< \epsilon_i$ 
small enough.
\end{prop}

We call the equivalence class of the unique configuration 
described by Proposition \ref{near_well-defined} (with
$0<\epsilon_i$ small enough) the {\it convex $l-$near-gon} or {\it near-gon}
with (cyclically oriented) near-edges $E_1,\dots,E_l$.
We denote by $P(E_1\cdots E_l)$ any configuration representing
the equivalence class of this near-gon. 

\begin{thm} \label{mainnear-gon} For every near-edge $E$ there exist 
polynomials
$$p(E)\in{\mathbf Z}[t]\hbox{ and }{\overline p}(E)\in{\mathbf Z}[s,t]$$
such that the convex near-polygon $P(E_1\cdots E_l)$ has
$$\big\langle \prod_{j=1}^l p(E_j),\sum_{n=2}^\infty C_{n-2}\ t^n\big\rangle_t
$$ 
maximal triangulations and complete triangulation polynomial
$$\big\langle \prod_{j=1}^l {\overline p}(E_j),
\sum_{n=2}^\infty C_{n-2}\ t^n\big\rangle_t$$
where $C_n={2n\choose n}/(n+1)$ are the Catalan numbers. 
\end{thm}

\begin{cor} The number of maximal triangulations and the complete
triangulation polynomial of a convex near-polygon $P(E_1,\dots,E_l)$
depend only on the multiset $\{E_1,\dots,E_l\}$ of (equivalence classes of)
its near-edges.
\end{cor}

\noindent{\bf Example.} The configuration of Figure $2$ is equivalent to 
the convex near-gon $P(E_a\cdot E_b\cdot E_c)$ where
$$\begin{array}{l}
\displaystyle
E_a=\left(\begin{array}{rrrrrr}0&1&2&3&4&5\\0&1&-1&1&-1&0\end{array}\right)\\
\displaystyle 
E_b=\left(\begin{array}{rrrrr}0&1&2&3&4\\0&1&-1&1&0\end{array}\right)\\
\displaystyle
E_c=\left(\begin{array}{rrrrrr}0&1&2&3&4&5\\0&2&1&-1&1&0\end{array}\right)
\end{array}$$
are near-edges whose complete edge-polynomials
$$\begin{array}{lcl}
\displaystyle
{\overline p}(E_a)&=&(14p_3+7p_4+p_5)s^5+(10p_2+7p_3+2p_4)s^4\\
\displaystyle &&\quad +(2p_1+2p_2+p_3)s^3\ ,\\
\displaystyle
{\overline p}(E_b)&=&(5p_3+p_4)s^4+2(2p_2+p_3)s^3+(p_1+p_2)s^2\\
\displaystyle
{\overline p}(E_c)&=&
(10p_3+7p_4+2p_5)s^5+(3p_2+13p_3+4p_4)s^4\\
\displaystyle &&\quad +3(2p_2+p_3)s^3+(p_1+p_2)s^2
\end{array}$$
(where $p_n=\sum_{k=0}^{\lfloor n/2\rfloor}(-1)^k
{n-k\choose k}t^{n-k}\in {
\mathbf Z}[t]$ as usual) will be computed in the sequel.
The complete triangulation polynomial
${\overline p}(E_a\cdot E_b\cdot E_c)$ of this configuration
is given by $$\begin{array}{cl}
\displaystyle &\langle {\overline p}(E_a){\overline p}(E_b)
{\overline p}(E_c),\sum_{n=2}^\infty C_{n-2}\ t^n\rangle_t\\
\displaystyle =&
194939\,{s}^{14}+338669\,{s}^{13}+263615\,{s}^{12}+119944\,{s}^{11}\\
\displaystyle 
&\quad +34773\,{s}^{10}+6522\,{s}^{9}+748\,{s}^{8}+42\,{s}^{7}\ .
\end{array}$$

\subsection{Edge-polynomials}

Given an $n-$near-edge 
$$E=(P_0,\dots,P_n)=\left(\begin{array}{cccc}
x_0&x_1&\dots&x_n\\ y_0&y_1&\dots&y_n\end{array}\right)$$
we define the {\it lower convex hull}  $\partial_- E$ as 
the piece-wise linear path joining $P_0$
to $P_n$ formed by the \lq\lq lower'' edges of $\partial \hbox{Conv}
(E)$. The set $V_-(E)$ of {\it lower extremal vertices} of $E$
is defined as the increasing subsequence $\hbox{Extr}(E)\cap \partial_- E$
of all extremal vertices of $\hbox{Conv}(E)$ situated in the closed 
\lq\lq lower'' halfplane containing $P_0$ and $P_n$ in its boundary
and containing all points $(x_0,y)$ with $y\leq y_0$.

A {\it roof} of $E$ is a strictly
increasing sub-sequence 
$$R=(P_{j_0}=P_0,P_{j_1},\dots,P_{j_{k-1}},P_{j_k}=P_n)\subset E$$
starting at $P_0$ and ending at $P_n$ of $E$. 
If a roof contains $k+1$ elements we call $\hbox{length}(R)=k$ 
its {\it length}.
The {\it skyline} $\hbox{Skyline}(R)$ of a roof $R$
is the piecewise-linear path from $P_0$ to $P_n$
defined by joining consecutive elements of $R$ using straight segments.

A roof $R$ is {\it covering} if it \lq\lq shelters'' every element of $E$:
for every point  $P_i=\left(\begin{array}{c} x_i\\y_i\end{array}\right)
\in E$ we have either $P_i\in R$ or 
$\tilde y_i> y_i$ for
$\left(\begin{array}{c} x_i\\{\tilde y}_i\end{array}\right)\in\hbox{Skyline}
(R)$.

The {\it maximal edge-polynomial} $p(E)$ of the near-edge $E$ is defined as
$$p(E)=\sum_{R\hbox{ covering roof of }E} \tau_{\hbox{max}}(E,R)\ 
p_{\hbox{length}(R)}\in{\mathbf Z}[t]$$
where $\tau_{\hbox{max}}(E,R)$ denotes the number of maximal 
triangulations of the (generally non-convex) compact polygonal region 
delimited by the two piecewise linear paths $\hbox{Skyline}(R)$ 
and $\partial_- E$.

A {\it $k-$sub-edge} of $E$ is an increasing subsequence 
$E'\subset E$ of $k+1\leq n+1$ elements in $E$ such that 
$V_-(E')=V_-(E)$. In particular, any sub-edge $E'$ of $E$ has also initial 
vertex $P_0$ and final vertex $P_n$. 

\begin{figure}[h]\label{8near-edge}
\epsfysize=1.8cm
\epsfbox{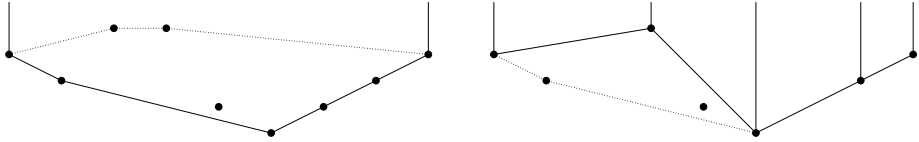}
\caption{An $8-$near-edge and a covering roof of weight $4$ of a
$6-$sub-edge}
\end{figure}  

\noindent{\bf Example.} The left side of Figure 7
%\ref{8near-edge}
displays the $8-$near edge
$$E=(P_0,\dots,P_8)=\left(\begin{array}{rrrrrrrrr}
0&1&2&3&4&5&6&7&8\\ 
0&-1&1&1&-2&-3&-2&-1&0\end{array}\right)\ .$$
We have $V_-(E)=(P_0,P_1,P_5,P_8)$
and $E$ has $2^5$ sub-edges obtained by removing any subset of 
vertices among $\{P_2,P_3,P_4,P_6,P_7\}$
from $E$. The right side of Figure 7
%\ref{8near-edge} 
displays
the covering roof $R=(P_0,P_3,P_5,P_7,P_8)$
of the sub-edge 
$E'=(P_0,P_1,P_3,P_4,P_5,P_7,P_8)\subset E$.

The {\it complete edge-polynomial} ${\overline p}(E)$ of an $n-$near-edge
$E=(P_0,\dots,P_n)$ is 
defined as
$${\overline p}(E)=\sum_{m=1}^n\sum_{E'\subset E\ m-\hbox{sub-edge}} p(E')s^m
$$
where $p(E')$ denotes the maximal edge-polynomial of the sub-edge
$E'\subset E$.

{\bf Example.} We compute the complete edge-polynomial 
${\overline p}(E_a)$ of the near-edge 
$$E_a=(P_0,\dots,P_5)=
\left(\begin{array}{rrrrrr}0&1&2&3&4&5\\0&1&-1&1&-1&0\end{array}\right)
$$
involved in the near-gon $P(E_a\cdot E_b\cdot E_c)$ of Figure 2.

\begin{figure}[h]\label{E5a}
\epsfysize=5.0cm
\epsfbox{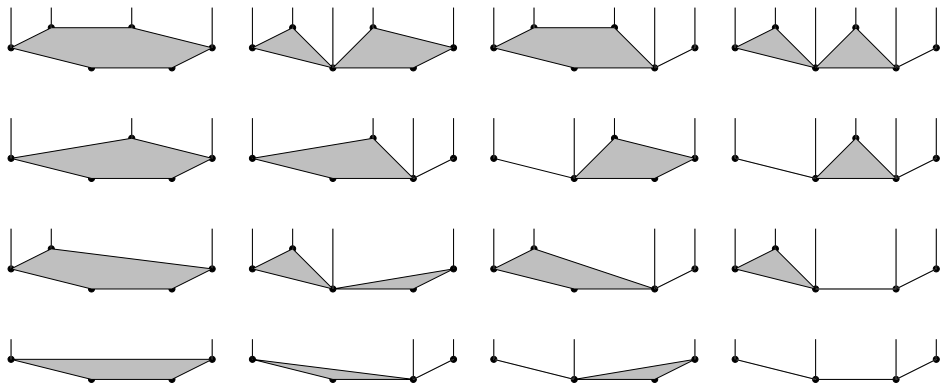}
\caption{All covering roofs of sub-edges involved in ${\overline p}(E_{a})$}
\end{figure}  

Figure 8
%\ref{E5a} 
contains all roofs of the four possible 
sub-edges of $E_a$ obtained by removing any subset of points in
$\{P_1,P_3\}$. Each possible sub-edge has
four different covering
roofs whose contributions to ${\overline p}(E_a)$ are given by
$$\begin{array}{|l||r|r|r|r|}
\hline
\hbox{sub-edge}& &&&\\
\hline\hline
\displaystyle (P_0,P_1,P_2,P_3,P_4,P_5):&14p_3s^5&2p_4s^5&5p_4s^5&p_5s^5\\
\hline\displaystyle (P_0,P_2,P_3,P_4,P_5):&5p_2s^4&2p_3s^4&2p_3s^4&p_4s^4\\
\hline\displaystyle (P_0,P_1,P_2,P_4,P_5):&5p_2s^4&p_3s^4&2p_3s^4&p_4s^4\\
\hline\displaystyle (P_0,P_2,P_4,P_5):&2p_1s^3&p_2s^3&p_2s^3&p_3s^3
\\ \hline\end{array}$$
They sum up to the complete edge-polynomial 
$${\overline p}(E_a)=(14p_3+7p_4+p_5)s^5+(10p_2+7p_3+2p_4)s^4+
(2p_1+2p_2+p_3)s^3$$ 
of $E_a$.

\subsection{Proof of Theorem \ref{mainnear-gon}}
Consider a triangulation ${\mathcal T}$ of a convex near-polygon
$P(E_1\cdots E_l)$ with complete triangulation polynomial
$${\overline p}=\sum_k\tau_k(P(E_1\cdots E_l))\ s^k\ .$$ 
A triangle $\Delta\in{\mathcal T}$ with vertices
$\{a,b,c\}$ is of {\it edge-type} if there exists $j$ such that
$\{a,b,c\}\subset E_j$. Otherwise, the triangle $\Delta$ is {\it 
interior}. If $\Delta$ with vertices $\{a,b,c\}$ is of edge-type, we write 
$\Delta\in E_j$ if $\{a,b,c\}\subset E_j$.

The definition of a convex near-polygon
implies that the set of all triangles $\Delta\in E_j\cap {\mathcal T}$ 
determines a covering roof $R_j$ of a suitable $k_j-$sub-edge $E'_j$ of
$E_j$. More precisely, define $E'_j$ as the set of all $(k_j+1)$ points of
$E_j$ which appear as vertices in ${\mathcal T}$. The set of edge-type
triangles in $E_j$ triangulates now a (generally non-convex) 
polygonal region enclosed by two (non-crossing) paths joining the
endpoints of $E_j$. These paths are the lower convex hull 
$\partial_-E'_j=\partial_-E_j$
and the skyline $\hbox{Skyline}(R_j)$ of a unique covering 
roof $R_j\subset E'_j$ having length $\hbox{length}(R_j)$.
The interior triangles of ${\mathcal T}$ define (after straightening)
a unique maximal triangulation of the weighted convex polygon
$P(\hbox{length}(R_1)\cdots \hbox{length}(R_l))$. The left side of Figure 9
%\ref{edgetriangle} 
displays a
maximal triangulation of a near-gon having four near-edges (of length
$2,3,1$ and $4$). The inner
triangles of the displayed triangulation yield a maximal triangulation
of the convex weighted polygon $P(2\cdot 2\cdot 1\cdot 3)$ depicted
on the right side of Figure 9.
%\ref{edgetriangle}.

\begin{figure}[h]\label{edgetriangle}
\epsfysize=3.2cm
\epsfbox{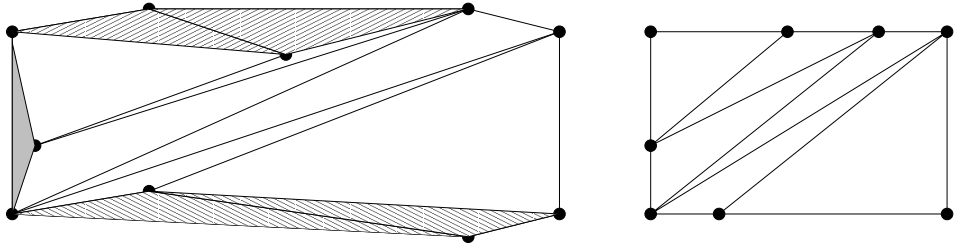}
\caption{Edge-type and interior triangles of a triangulation}
\end{figure}  

There are exactly 
$$\tau_{\hbox{max}}(\hbox{length}(R_1)\cdots \hbox{length}(R_l))\ 
\prod_{j=1}^l\tau_{\hbox{max}} (E'_j,R_j)$$
such triangulations contributing to the coefficient
$s^{(k_1+k_2+\dots+k_l)}$ of 
$\overline p$ for fixed covering roofs $R_j$ of fixed 
$k_j-$sub-edges $E'_j\subset E_j$. Applying Theorem \ref{maxformulaweighted} 
and summing over all covering roofs in sub-edges
yields the result.\hfill $\Box$

%%%%%%%%%%%%%%%%%%%%%%%%%%%%%%%%%%%%%%%%%%%%%%%%%%%%%%%%%%%%%%%%%%%%%%%%%

\section{Properties of edge-polynomials}

A near edge $E=(P_0,\dots,P_n)$ {\it factorises} into near-edges
$E_1,E_2$ if there exists a lower extremal vertex
$P_k\in E\cap \partial_- E$ such that 
$$E_1=(P_0,P_1,\dots,P_{k-1},P_k),\ E_2=(P_k,P_{k+1},\dots,P_{n-1},P_n)$$
and all points of $E\setminus E_i$ lie strictly above any line defined 
by two distinct points of $E_i$ for $i=1,2$. We write $E=E_1\cdot E_2$
if the near-edge $E$ factorises with first factor $E_1$ and second 
factor $E_2$. A near-edge is {\it prime} if it has no non-trivial
factorisation. It is easy to show that every near-edge 
has a unique factorisation into prime near-edges.

\begin{prop} Given a factorisation $E=E_1\cdot E_2$ of a near-edge $E$
we have
$$p(E)=p(E_1)\ p(E_2)\hbox{ and }{\overline p}(E)={\overline p}(E_1)\ 
{\overline p}(E_2)\ .$$
\end{prop}

\noindent{\bf Proof.} Since the 
near-polygons $P(1^k\cdot E)$ and $P(1^k\cdot E_1\cdot E_2)$
are equivalent for all $k=2,3,\dots$ we have 
$$\big\langle t^k\ {\overline p}(E),\sum_{n=2}^\infty
C_{n-2}\ t^n\big\rangle_t=\big\langle t^k\ {\overline p}(E_1)\ 
{\overline p}(E_2),\sum_{n=2}^\infty
C_{n-2}\ t^n\big\rangle_t\ .$$
This implies the result
since $\det(\big(C_{i+j+k}\big)_{0\leq i,j\leq n})>0$
for all $k\geq 0$ and $n\geq 1$ (this follows for instance easily from 
Exercice 6.26.b in \cite{St}).\hfill $\Box$

Since the so-called $k-$shifted Hankel matrices 
$$H_{i,j}=C_{i+j+k},\ 0\leq i,j<n$$
are all non-singular, 
the coefficients $\gamma_\alpha,\dots,\gamma_d$ (with $\alpha\geq 2$)
of a polynomial $q(t)=\sum_{i=\alpha}^d \gamma_i\ t^i$ can be obtained 
from the $(d-\alpha+1)$ \lq\lq linear form evaluations''
$$\langle t^k\ q(t),\sum_{n=2}^{\infty}C_{n-2}\ t^n\rangle_t,\ 
k=0\dots,d-\alpha\ .$$

This turns any method or algorithm for enumerating 
maximal (or all) triangulations of configurations into 
a method or algorithm for computing maximal (or complete)
edge-polynomials.

\noindent{\bf Example.}
Using the online computing service \cite{Aol} of O. Aichholzer,
we compute the
maximal triangulation polynomial $p(E_c)$ (which is of of course
also the coefficient of $s^5$ in the complete triangulation
polynomial ${\overline p}(E_c)$) of the near-edge
$$E_c=\left(\begin{array}{rrrrrr}0&1&2&3&4&5\\0&2&1&-1&1&0\end{array}
\right)$$
involved in the near-gon of Figure $2$.

Since the set $\partial_+ (E_c)$ of all boundary vertices
on the upper boundary of $\hbox{Conv}(E_c)$ 
consists of the $3$ segments $\overline{P_0 P_1},\ \overline{P_1P_4}$
and $\overline{P_4P_5}$, we have
$$p(E_c)=\alpha p_3+\beta p_4+\gamma p_5$$
where $\alpha,\beta$ and $\gamma$ are unknown natural integers.

We identify now the near-edge $E_c=(P_0,\dots,P_5)$ with its underlying set
$$\{\left(\begin{array}{r}0\\0\end{array}\right),
\left(\begin{array}{r}1\\2\end{array}\right),
\left(\begin{array}{r}2\\1\end{array}\right),
\left(\begin{array}{r}3\\-1\end{array}\right),
\left(\begin{array}{r}4\\1\end{array}\right),
\left(\begin{array}{r}5\\0\end{array}\right)\}\ .$$
Define
$$A_1=\left(\begin{array}{r}1\\10\end{array}\right),\ 
A_2\left(\begin{array}{r}2\\11\end{array}\right),
A_3\left(\begin{array}{r}3\\10\end{array}\right)$$
and check that the configurations
$$E_c\cup\{ A_1\},\ E_c\cup\{A_1,A_2\},\ E_c\cup\{A_1,A_2,A_3\}$$
are near-polygons with near-edge factorisation
$$E_c\cdot E_{1}^2,E_c\cdot E_{1}^3,E_c\cdot E_{1}^4$$
where $E_{1}$ denotes the unique $1-$edge with maximal edge polynomial 
$p(E_1)=p_1=t$.
Denoting by $\tau_{\hbox{max}}({\mathcal C})$ the number of maximal 
triangulations of a finite configuration ${\mathcal C}$ and using for 
instance the on-line computing service of Aichholzer \cite{Aol}
mentionned above we get the linear system of equations
$$\begin{array}{rclcl}
19&=&\tau_{\hbox{max}}(E_c\cup\{ A_1\})&=&\langle p(E_c)\ t^2,
\sum_{n=2}^\infty C_{n-2}\ t^n\rangle_t\\
87&=&\tau_{\hbox{max}}(E_c\cup\{ A_1,A_2\})&=&\langle p(E_c)\ t^3,
\sum_{n=2}^\infty C_{n-2}\ t^n\rangle_t\\
175&=&\tau_{\hbox{max}}(E_c\cup\{ A_1,A_2,A_3\})&=&\langle p(E_c)\ t^4,
\sum_{n=2}^\infty C_{n-2}\ t^n\rangle_t\\ \end{array}$$
Substituting $p_n=
\sum_{k=0}^{\lfloor n/2\rfloor}(-1)^k{n-k\choose k}t^{n-k},\
C_n={2n\choose n}/(n+1)$ and solving for the unknowns
$\alpha,\beta,\gamma$ we get $$p(E_c)=10p_3+7p_4+2p_5\ .$$

\noindent{\bf Remark.} Since the computations of the coefficient
$\alpha=\tau_{\hbox{max}}(E_c)=10$ of lowest order (in $t$) and
of the leading coefficient  
$$
\gamma=
\tau_{\hbox{max}}(\{
\left(\begin{array}{r}0\\0\end{array}\right),
\left(\begin{array}{r}1\\2\end{array}\right),
\left(\begin{array}{r}2\\1\end{array}\right),
\left(\begin{array}{r}3\\-1\end{array}\right)\})\ \cdot $$
$$
\qquad \cdot\ 
\tau_{\hbox{max}}(\{\left(\begin{array}{r}3\\-1\end{array}\right),
\left(\begin{array}{r}4\\1\end{array}\right),
\left(\begin{array}{r}5\\0\end{array}\right)
\})=2\cdot 1=2$$
are easy, we could have simplified the above 
computation by a considerable amount of work.

For completeness we mention also the follwing obvious fact:

Call two $n-$near-edges $E=(P_0,\dots,P_n)$ and $E'=(P'_0,\dots,p'_n)$
{\it vertical mirrors} if $E'$ given by
$E'=(\overline P_n,\overline P_{n-1},\dots,\overline P_1,\overline P_0)$
where $\overline P_i=\left(\begin{array}{c}-x_i\\y_i\end{array}\right)$
is the Euclidean reflection of 
$P_i=\left(\begin{array}{c}x_i\\y_i\end{array}\right)$ with respect 
to the vertical line $x=0$. 

The following result is obvious:

\begin{prop} If a pair of near-edges $(E,E')$ are vertical mirrors,
then
$$\overline p(E)=\overline p(E')\ .$$
\end{prop}

\subsection{Complete edge-polynomials for small near-edges}

This subsection contains representants and edge-polynomials
for all $1-$, $2-$ and $3-$near-edges.

\subsubsection{$1-$near-edges} The unique $1-$near-edge can be
represented by $E_{1}=\left(\begin{array}{cc}0&1\\0&0\end{array}\right)$.
It is generic and prime and
has complete edge-polynomial ${\overline p}(E_{1})=p_1s=s\ t$.

\subsubsection{$2-$near-edges} There are two generic $2-$near-edges,
represented by
$$E_{2,1}=\left(\begin{array}{rrr}0&1&2\\0&1&0\end{array}\right),\quad
E_{2,2}=\left(\begin{array}{rrr}0&1&2\\0&-1&0\end{array}\right)\ .$$
$E_{2,1}$ is prime while $E_{2,2}=E_{1}\cdot E_{1}$. 
They have complete edge-polynomials
$${\overline p}(E_{2,1})=p_2s^2+p_1s,\ {\overline p}(E_{2,2})=
(p_2+p_1)s^2={\overline p}(E_{1})^2\ .$$

Moreover, there is also a unique non-generic $2-$near-edge represented 
for instance by 
$$E_{2,3}=\left(\begin{array}{rrr}0&1&2\\0&0&0\end{array}\right)$$
with complete edge-polynomial given by
$${\overline p}(E_{2,3})=p_2s^2+p_1s={\overline p}(E_{2,1})\ .$$

\subsubsection{$3-$near-edges} There are (up to equivalence) eight
generic $3-$near edges represented by
$$\begin{array}{ll}
\displaystyle 
E_{3,1}=\left(\begin{array}{rrrr}0&1&2&3\\0&1&3&0\end{array}\right)\quad ,&
E_{3,2}=\left(\begin{array}{rrrr}0&1&2&3\\0&1&1&0\end{array}\right)\ ,\\
\displaystyle 
E_{3,3}=\left(\begin{array}{rrrr}0&1&2&3\\0&3&1&0\end{array}\right)\ ,& 
E_{3,4}=\left(\begin{array}{rrrr}0&1&2&3\\0&1&-1&0\end{array}\right)\quad ,\\
\displaystyle 
E_{3,5}=\left(\begin{array}{rrrr}0&1&2&3\\0&-1&1&0\end{array}\right)\ ,&
E_{3,6}=\left(\begin{array}{rrrr}0&1&2&3\\0&-3&-1&0\end{array}\right)\quad ,\\
\displaystyle
E_{3,7}=\left(\begin{array}{rrrr}0&1&2&3\\0&-1&-1&0\end{array}\right)\quad ,&
E_{3,8}=\left(\begin{array}{rrrr}0&1&2&3\\0&-1&-3&0\end{array}\right)\ 
\quad .\\
\end{array} 
$$

The first five are prime. The last three have factorisations 
$$E_{3,6}=E_{1}E_{2,1}\ ,E_{3,7}=E_{1}^3,\ E_{3,8}=E_{2,1}E_{1}\ .$$
The pairs $\{E_{3,1},E_{3,3}\},\ \{E_{3,4},E_{3,5}\},\ \{E_{3,6},E_{3,8}\}$
are vertical mirrors. The prime near-edges have complete polynomials
$$\begin{array}{l}
\displaystyle {\overline p}(E_{3,1})={\overline p}(E_{3,3})=
(p_2+p_3)s^3+2p_2s^2+p_1s\ ,\\
\displaystyle {\overline p}(E_{3,2})=2p_3s^3+2p_2s^2+p_1s\ ,\\ 
\displaystyle {\overline p}(E_{3,4})={\overline p}(E_{3,5})=
(2p_2+p_3)s^3+p_1^2s^2\ .
\end{array}$$

There are moreover nine more $3-$near-edges which are not generic. They 
are represented for instance by

$$\begin{array}{ll}
\displaystyle 
E_{3,9}=\left(\begin{array}{rrrr}0&1&2&3\\0&1&2&0\end{array}\right)\quad ,&
E_{3,10}=\left(\begin{array}{rrrr}0&1&2&3\\0&-1&-2&0\end{array}\right)\ ,\\
\displaystyle 
E_{3,11}=\left(\begin{array}{rrrr}0&1&2&3\\0&0&1&0\end{array}\right)\ ,& 
E_{3,12}=\left(\begin{array}{rrrr}0&1&2&3\\0&0&-1&0\end{array}\right)\quad ,\\
\displaystyle 
E_{3,13}=\left(\begin{array}{rrrr}0&1&2&3\\0&1&0&0\end{array}\right)\ ,&
E_{3,14}=\left(\begin{array}{rrrr}0&1&2&3\\0&-1&0&0\end{array}\right)\quad ,\\
\displaystyle
E_{3,15}=\left(\begin{array}{rrrr}0&1&2&3\\0&2&1&0\end{array}\right)\quad ,&
E_{3,16}=\left(\begin{array}{rrrr}0&1&2&3\\0&-2&-1&0\end{array}\right)\ 
\quad ,\\
\displaystyle
E_{3,17}=\left(\begin{array}{rrrr}0&1&2&3\\0&0&0&0\end{array}\right)
\end{array} 
$$

The near-edges 
$$E_{3,10}=E_{2,3}\ E_{1},\qquad E_{3,16}=E_{1}\ E_{2,3}$$
factorise. The remaining prime near-edges have complete
edge-polynomials

$$\begin{array}{lclcl}
\displaystyle 
{\overline p}(E_{3,9})&=&{\overline p}(E_{3,15})&=&p_3s^3+2p_2s^2+p_1s,\\
\displaystyle
{\overline p}(E_{3,11})&=&{\overline p}(E_{3,13})&=&
(p_2+p_3)s^3+2p_2s^2+p_1s,\\
\displaystyle 
{\overline p}(E_{3,12})&=&{\overline p}(E_{3,14})&=&(p_2+p_3)s^3+p_2s^2,\\
\displaystyle 
{\overline p}(E_{3,17})&=&{\overline p}_3&=&p_3s^3+2p_2s^2+p_1s\ .
\end{array}$$

\subsection{Enumeration of generic near-edges}

This subsection is a digression sketching briefly a bijection
between near-edges (up to equivalence) and orbits of
suitable extremal vertices of configurations under automorphisms.  

Let $P\in\hbox{Extr}({\mathcal C})$ be an extremal vertex of a configuration
${\mathcal C}=\{P_1,\dots,P_{n+2}\}\subset {\mathbf R}^2$
consisting of $(n+2)$ points. Suppose that $P$ is
not collinear with two distinct elements of ${\mathcal C}\setminus
\{P\}$. Let $P_-,P_+\in\hbox{Extr}({\mathcal C})$ be the 
immediate predecessor and successor (in counterclockwise order) of $P$
on the boundary $\partial\hbox{Conv}({\mathcal C})$.

A suitable projective transformation (sending a line outside
$\hbox{Conv}({\mathcal C})$ which misses $P$ only nearly to the
line at infinity and sending $P$ close to the ideal point $(0,+\infty)$) 
transforms ${\mathcal C}\setminus\{P\}$ into an $n-$near-edge
$E$ with initial vertex $P_+$ and final vertex $P_-$.  
This near-edge $E$ is well-defined, up to equivalence, and 
every $n-$near-edge can be obtained in this way.

\begin{figure}[h]\label{3near_edges}
\epsfysize=3.6cm
\epsfbox{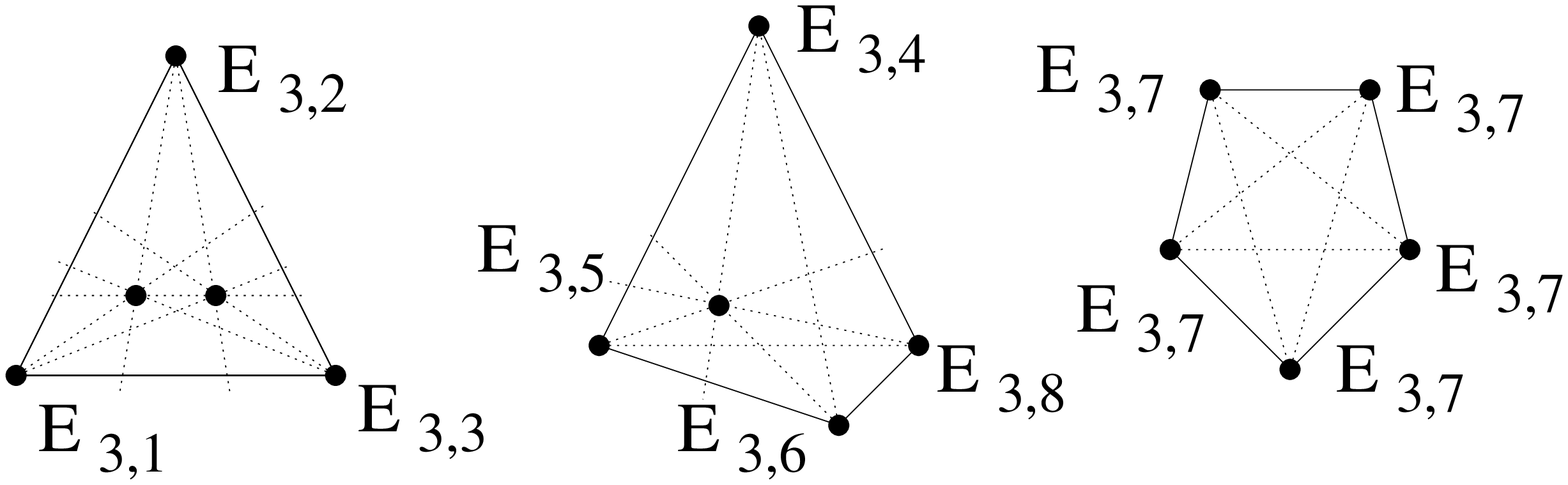}
\caption{Generic configurations of $5$ points and the associated generic 
$3-$near-edges}
\end{figure}  

In particular, generic $n-$near-edges are in bijection with orbits under
automorphisms of extremal vertices in generic configurations of
$(n+2)$ points. Figure 10
%\label{3near_edges} 
illustrates this
for $n=3$. It shows all three generic configurations of five
points with extremal vertices labelled by the corresponding generic
$3-$near-edges (using the notations of the previous subsection).

Call two generic $n-$near-edges {\it related} if they are associated
to two extremal vertices of the same generic configuration of $(n+2)$
points. Classes of related generic $n-$near-edges are in bijection
with equivalence classes of generic configurations having $(n+2)$ points.

We say that two configurations $E=(P_0,\dots,P_n)$ and
$E'=(P'_0,\dots,P'_n)$ are {\it related by a horizontal reflection} 
if $P'_i=\left(\begin{array}{c}x_i\\ -y_i\end{array}\right)$ is
the reflection with respect to the horizontal line $y=0$ of 
$P_i=\left(\begin{array}{c}x_i\\ y_i\end{array}\right)$.

Equivalence classes of related generic $n-$near-edges,
up to horizontal reflections, are in bijection with generic projective 
configurations (up to equivalence) of $(n+2)$ points in the
projective plane ${\mathbf R}P^2$.

%%%%%%%%%%%%%%%%%%%%%%%%%%%%%%%%%%%%%%%%%%%%%%%%%%%%%%%%%%%%%%%%%%%%%%%%%%%%%

\section{An algorithm for convex near-edges}

An $n-$near-edge $E=\{P_0,\dots,P_n\}$ 
is {\it convex} if its points $P_0,\dots,P_n$ are all extremal, i.e.
vertices of a convex polygon with $(n+1)$ edges. 
This section describes an algorithm for computing 
maximal and complete edge-polynomials of convex near-edges.

A convex $n-$near-edge can be represented by the sequence of points
$$(P_0=\left(\begin{array}{c}0\\ 0\end{array}\right),\dots,
P_i=\left(\begin{array}{c}i\\ \epsilon_i\ i(n-i)\end{array}\right),\dots,
P_n=\left(\begin{array}{c}n\\ 0\end{array}\right))$$
where $\epsilon_1,\dots,\epsilon_{n-1}\in\{\pm 1\}$. There are
$2^{n-1}$ equivalence classes,
encoded by $(\epsilon_1,\dots,\epsilon_{n-1})\in\{\pm 1\}^{n-1}$, 
of convex $n-$near-edges. The 
convex near-edge with $\epsilon_i=-1,\ i=1,\dots,(n-1)$ has the factorisation
$E_{1}^n$. All others are prime.

The following algorithm of polynomial complexity
computes the complete (respectively maximal)
edge-polynomial of a convex near-edge. It is based on the fact that 
all connected components of the open polygonal region delimited by
a covering roof $R$ and $\partial_-(E)$ are open convex polygons if $E$ is a 
convex near-edge. The exponent of $u$ keeps track of the
length of a \lq\lq covering roof under construction'' 
$(P_0,P_{i_1},\dots,P_{i_j},?)$. The exponent of $w$ keeps
track of the number of vertices of the rightmost polygonal region
which has not yet closed up. 

\noindent{\bf Algorithm} 

Input: $n\in 1,2,\dots,$ and $(\epsilon_1,\epsilon_2,\dots,
\epsilon_{n-1})\in\{\pm 1\}^{n-1}$.

Set 
$$r= \left\{\begin{array}{ll}
\displaystyle s\cdot u\qquad &\hbox{(complete edge-polynomial)}\\
\displaystyle u&\hbox{(maximal edge-polynomial)}
\end{array}\right.$$

For $i$ from 1 to $n-1$ do $\{$

If $\epsilon_i=1$ replace

$$r\longmapsto \left\{\begin{array}{ll}
\displaystyle r\cdot (s\cdot u\cdot w+1)&\hbox{(complete edge-polynomial)}\\
\displaystyle r\cdot u\cdot w&\hbox{(maximal edge-polynomial)}
\end{array}\right.$$

If $\epsilon_i=-1$ replace

$$r\longmapsto \left\{\begin{array}{ll}
\displaystyle r\cdot s\cdot w+\big\langle\ r,\left(\sum_{n=0}^\infty C_n\ w^n
\right)\big\rangle_w\cdot s\cdot u&\hbox{(complete edge-polynomial)}\\
\displaystyle r\cdot w+\langle r,\left(\sum_{n=0}^\infty C_n\ w^n
\right)\rangle_w\ u
&\hbox{(maximal edge-polynomial)}
\end{array}\right.$$

$\}$

Result:

$$\big\langle r,\left(\sum_{n=0}^\infty C_n\ w^n\right)\ \left(
\sum_{j=1}^\infty p_j\ u^j\right)\big\rangle_{u,w}$$
(where $p_n=\sum_k (-1)^k{n-k\choose k}t^{n-k}$).

%\begin{prop} The result of the above algorithm is the complete
%or maximal edge-polynomial of the convex $n-$near-edge
%encoded by $(\epsilon_1,\dots,\epsilon_{n-1})$.
%\end{prop}

\noindent{\bf Example}. The first two near-edges 
$$
E_a=\left(\begin{array}{rrrrrr}0&1&2&3&4&5\\0&1&-1&1&-1&0\end{array}\right)$$
and
$$
E_b=\left(\begin{array}{rrrrr}0&1&2&3&4\\0&1&-1&1&0\end{array}\right)$$
involved in Figure $2$ are convex. 
The near-edge $E_a$ corresponds to $n=5$ and the sequence 
$(1,-1,1,-1)$. The near-edge $E_b$ corresponds to $n=4$ and $(1,-1,1)$.

The above algorithm for $E_a$ 
(with $n=5$ and $(\epsilon_1,\dots,\epsilon_4)=(1,-1,1,-1)$) yields

$$
\begin{array}{|c|r|l|}
\hline
i=&\epsilon_i=&r=\\
\hline\hline
&&su\\
1&1&s^2u^2w+su\\
2&-1&s^3u^2w^2+s^2uw+s^3u^3+s^2u^2\\
3&1&s^4u^3w^3+s^3u^2w^2+s^4u^4w+s^3u^3w+\\
&&\quad +s^3u^2w^2+s^2uw+s^3u^3+s^2u^2\\
4&-1&
s^5u^3w^4+2s^4u^2w^3+s^5u^4w^2+s^4u^3w^2+s^3uw^2+s^4u^3w\\
&&\quad +s^3u^2w+5s^5u^4+4s^4u^3+s^5u^5+2s^4u^4+s^3u^2+s^3u^3\\
\hline\end{array}$$
$$\hbox{Result: }(14p_3+7p_4+p_5)s^5+(10p_2+7p_3+2p_4)s^4+(2p_1+2p_2+p_3)s^3$$
which is the complete edge-polynomial ${\overline p}(E_a)$ of $E_a$.

Computing 
$$\big\langle r,\left(\sum_{n=0}^\infty C_n\ w^n\right)\ \left(
\sum_{j=1}^\infty p_j\ u^j\right)\big\rangle_{u,w}$$
with
$$r= s^4u^3w^3+2s^3u^2w^2+s^4u^4w+s^3u^3w+s^2uw+s^3u^3+s^2u^2$$
corresponding to $i=3$
we get the complete edge-polynomial
$${\overline p}(E_b)=(5p_3+p_4)s^4+(4p_2+2p_3)s^3+(p_1+p_2)s^2$$
of the convex near-edge $E_b$ encoded by $n=4$ and $(1,-1,1)$.

%%%%%%%%%%%%%%%%%%%%%%%%%%%%%%%%%%%%%%%%%%%%%%%%%%%%%%%%%%%%%%%%%%%%%%%%%%%%

\section{Roof-sequences and a transfer matrix}

In this section, we define the set of roof-sequences of a finite
configuration (or near-edge). Such sequences are in
bijection with triangulations. Using this bijection, 
we introduce a transfer matrix for computing
the complete triangulation polynomial (or complete edge-polynomial)
of a configuration (or near-edge).  

We endow ${\mathbf R}^2$ with  the total order relation defined by 
$$\left(\begin{array}{c}x\\y\end{array}\right)<
\left(\begin{array}{c}x'\\y'\end{array}\right)$$
if either $x<x'$ or $x=x'$ and $y>y'$. This order is obtained by
reflecting (orthogonaly) the plane endowed with the
familiar lexicographic order
with respect to the line $y=0$.

Consider a finite configuration ${\mathcal C}=\{P_0,\dots,P_n\}
\subset {\mathbf R}^2$ consisting of $(n+1)$ totally ordered
points $P_0<P_1<\dots <P_n$ for the above order. We 
identify ${\mathcal C}$ with the strictly increasing sequence 
$(P_0,\dots , P_{n})$
of its $(n+1)$ points.
The initial point $P_0$ and the final point $P_n$ of 
${\mathcal C}$ are extremal vertices of $\hbox{Conv}({\mathcal C})$
and are joined
by two (generally distinct) piecewise-linear paths $\partial_+({\mathcal C})$
and $\partial_-({\mathcal C})$, called the {\it upper} and {\it lower}
convex hull. The path $\partial_+({\mathcal C})$ (respectively
$\partial_-({\mathcal C})$) is contained in the
boundary $\partial \hbox{Conv}({\mathcal C})$ and links $ P_0$ to $P_n$ 
by turning clockwise (respectively counterclockwise) 
around the convex hull 
$\hbox{Conv}({\mathcal C})$ of ${\mathcal C}$.

A {\it roof} $R$ of ${\mathcal C}$ having length $w=\hbox{length}(R)$ 
is an increasing subsequence
$R=(P_{i_0}=P_0,P_{i_1},\dots,P_{i_{w-1}},P_{i_w}=P_n)$ containing
$(w+1)$ elements.

The {\it skyline} of a (decorated) roof $R=(P_0,\dots,P_n)$ is the
piece-wise linear path from $P_0$ to $P_n$ obtained by
joining consecutive points of $R$ using rectilinear segments.
Skylines and roofs are in 
bijection for generic configurations. If ${\mathcal C}$ is non-generic,
several roofs might share a common skyline. 

\begin{figure}[h]\label{expleboundandroof}
\epsfysize=5cm
\epsfbox{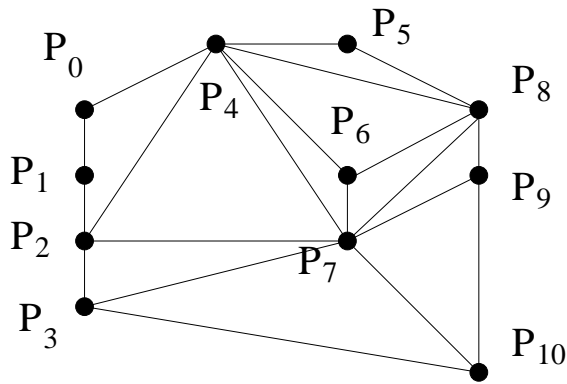}
\caption{A triangulated configuration of $(10+1)$ points}
\end{figure}  

\noindent{\bf Example.} The configuration 
$${\mathcal C}=\{P_0,\dots,P_n\}=\left(\begin{array}{rrrrrrrrrrrrrrrr}
0&0&0&0&1&2&2&2&3&3&3\\
3&2&1&0&4&4&2&1&3&2&-1\end{array}\right)$$
represented in Figure 11 has four roofs with skyline 
$\partial_-{\mathcal C}$:
$$\begin{array}{ll}
\displaystyle (P_0,P_3,P_{10})\ ,&(P_0,P_1,P_3,P_{10})\ ,\\
\displaystyle (P_0,P_2,P_3,P_{10})\ ,\quad &(P_0,P_1,P_2,P_3,P_{10})
\end{array}$$
and two roofs with skyline $\partial_+{\mathcal C}$:
$$(P_0,P_4,P_5,P_8,P_{10})\,\quad (P_0,P_4,P_5,P_8,P_9,P_{10})\ .$$

A {\it decorated roof} is a roof 
$$R=(P_{i_0}=P_0,P_{i_1},\dots,P_{i_{w-1}},P_{i_w}=P_n)\ ,$$
decorated by one of the
segments joining two consecutive points $P_{i_d},P_{i_{d+1}}$ of $R$.
We denote such a decorated roof by 
$$R=(P_0,\dots,P_{i_{d-1}},\overline{P_{i_d},P_{i_{d+1}}},P_{i_{d+2}},\dots,
P_{i_w}=P_n)$$ and call $R$ {\it initially decorated} if $d=0$, i.e. if
the decoration starts at the initial point $P_0$. 

Three non-collinear points $A<B<C\in {\mathcal C}$ are vertices of 
a triangle {\it of type V}
if $B$ is \lq\lq above'' the line containing $A$ and $C$. The 
non-collinear points $A,B,C$ are
vertices of a triangle {\it of type $\Lambda$} otherwise. In other terms,   
a triangle with (necessarily non-collinear) vertices $A<B<C$ is of type
$\Lambda$ if $\det(C-A,B-A)>0$ and of type $V$ if $\det(C-A,B-A)<0$.

We call a triangle with non-collinear vertices $A,B,C$ {\it minimal} if 
$\{A,B,C\}={\mathcal C}\cap \hbox{Conv}(A,B,C)$. In particular,
a triangulation is maximal if and only if all its triangles are
minimal.

A decorated roof
$$R'=(P_{i_0}=P_0,\dots,P_{i_{k-1}},\overline{P_{i_k},P_{i_{k+1}}},
P_{i_{k+2}},P_{i_{k+3}},\dots,P_{i_w}=P_n)$$
of length $w$
is a {\it $\Lambda-$successor} of $R$ 
if $R$ is a decorated roof of length $w-1$ given by 
$$R=(P_{i_0}=P_0,\dots,\overline{P_{i_d},P_{i_{d+1}}},\dots,
P_{i_{k-1}},P_{i_k},P_{i_{k+2}},P_{i_{k+3}},\dots,P_{i_w}=P_n)$$
satisfying $P_{i_d}\leq P_{i_k}$ and the points $\{P_{i_k},P_{i_{k+1}},
P_{i_{k+2}}\}\subset R'$ are vertices of a $\Lambda-$triangle. 

Otherwise stated, $R'$ is obtained from $R$ by gluing a $\Lambda-$triangle
with decorated left-upper edge onto an edge of $R$ which is not situated
to the left of the decorated edge in $R$.

Similarly,  a decorated roof
$$R'=(P_{i_0}=P_0,\dots,P_{i_{k-1}},\overline{P_{i_k},P_{i_{k+2}}},
P_{i_{k+3}},\dots,P_{i_w}=P_n)$$
of length $w-1$
is a {\it $V-$successor} of $R$ if $R$ is the decorated 
roof of length $w$ given by 
$$R=(P_{i_0}=P_0,\dots,\overline{P_{i_d},P_{i_{d+1}}},\dots,
P_{i_k},P_{i_{k+1}},
P_{i_{k+2}},P_{i_{k+3}},\dots,P_{i_w}=P_n)$$
satisfying $P_{i_{d+1}}\leq P_{i_{k+2}}$ and the points $\{P_{i_k},P_{i_{k+1}},
P_{i_{k+2}}\}\subset R$ are vertices of a $V-$triangle. 

Otherwise stated, $R'$ is obtained from $R$ by gluing a $V-$triangle
with decorated upper edge onto two consecutive edges of increasing slope 
of $R$ and the decorated edge of $R$ is not situated to the right of these
two consecutive edges.

\begin{figure}[h]\label{LambdaVtriangles}
\epsfysize=10cm
\epsfbox{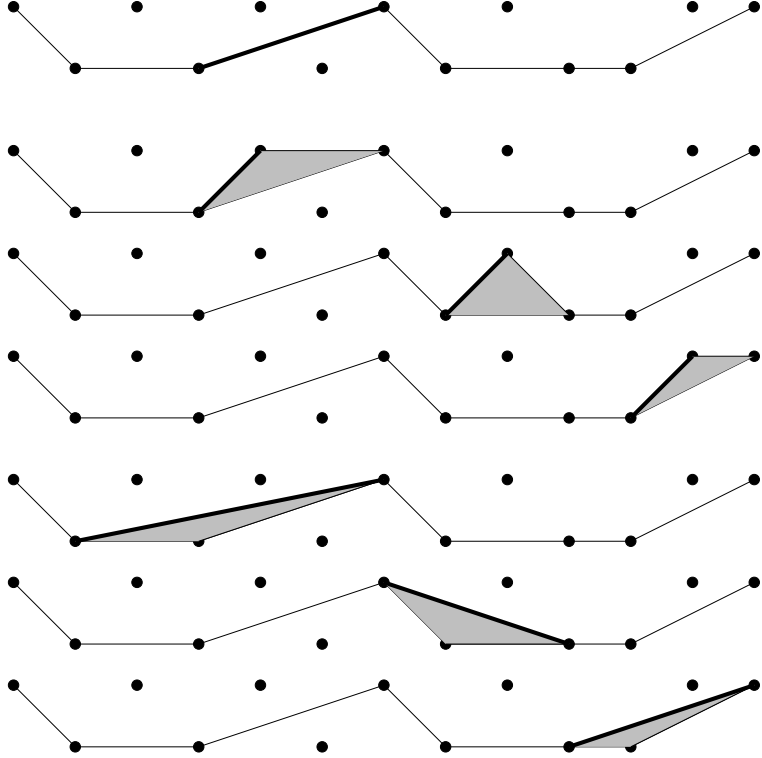}
\caption{A decorated roof and all its $\Lambda-$ and $V-$successors}
\end{figure}  

\noindent{\bf Example.} Figure 12
%\ref{LambdaVtriangles}
displays the decorated roof
$$R=(P_0,P_1,\overline{P_3,P_6},P_7,P_9,P_{10},P_{12})$$ 
of the $12-$near-edge
$$E=\left(\begin{array}{rrrrrrrrrrrrrrr}
0&1&2&3&4&5&6&7&8&9&10&11&12\\
1&0&1&0&1&0&1&0&1&0&0&1&1\end{array}\right)\ .$$
and all its $\Lambda-$ and $V-$successors
$$\begin{array}{l}
(P_0,P_1,\overline{P_3,P_4},P_6,P_7,P_9,P_{10},P_{12})\\
(P_0,P_1,P_3,P_6,\overline{P_7,P_8},P_9,P_{10},P_{12})\\
(P_0,P_1,P_3,P_6,P_7,P_9,\overline{P_{10},P_{11}},P_{12})\\
(P_0,\overline{P_1,P_6},P_7,P_9,P_{10},P_{12})\\
(P_0,P_1,P_3,\overline{P_6,P_9},P_{10},P_{12})\\
(P_0,P_1,P_3,P_6,P_7,\overline{P_9,P_{12}})
\end{array}$$

A decorated roof $R'$ is a {\it successor} of a 
decorated roof $R$ if $R'$ is either a 
$\Lambda-$ or a $V-$successor of $R$. We call $R'$ an {\it immediate
successor} of $R$ if the difference between $R$ and $R'$ is a minimal
triangle.

Given a finite configuration ${\mathcal C}\subset {\mathbf R}^2$, a 
{\it (maximal) roof-sequence} is a sequence 
$$R_{i_1},\dots,R_{i_l}$$
of decorated roofs where $R_{i_1}$ is initially decorated
and where $R_{i_{j+1}}$ is a (immediate) successor of 
$R_{i_j}$ for $j=1,\dots,l-1$. The skyline of $R_{i_1}$ (respectively
of $R_{i_l}$) is the {\it initial} (respectively {\it final})
skyline of the roof-sequence.

\begin{prop} \label{maxroofseq}
The set of all (maximal) roof-sequences $R_{i_1},\dots,R_{i_l}$
with initial skyline $S_\alpha$ and final skyline $S_\omega$ is in bijection
with the set of all (maximal)
triangulations with vertices in ${\mathcal C}$
of the polygonal region enclosed by $S_\alpha$ and $S_\omega$.
\end{prop}

\noindent{\bf Proof.} Given a roof-sequence 
$(R_{i_1},\dots,R_{i_l})$, the union of all associated
skylines triangulates (maximally) the polygonal region enclosed by 
$S_{i_1}$ and $S_{i_l}$.

Given a (maximal) triangulation of a polygonal region enclosed by
two skylines, remove vertically the rightmost triangle
which has only \lq\lq sky'' above it (and which is a minimal
triangle for a maximal triangulation). Iterating this
algorithm one gets an undecorated roof-sequence in reversed 
order. The decorations can now be added only in a unique
way, except for the lowest skyline which we decorate initially.
\hfill $\Box$

\noindent{\bf Example.} The triangulation of the configuration
shown in Figure 11 corresponds to the roof-sequence 
$$\begin{array}{l}
\displaystyle (\overline{P_0,P_2},P_3,P_{10}),\ (P_0,P_2,\overline{
P_3,P_7},P_{10}),\ (P_0,\overline{P_2,P_7},P_{10}),\\
\displaystyle (P_0,\overline{P_2,P_4},P_7,P_{10}),\ 
(\overline{P_0,P_4},P_7,P_{10}),\ (P_0\overline{P_4,P_6},P_7,P_{10}),\\
\displaystyle (P_0,P_4,P_6,\overline{P_7,P_9},P_{10}),\
(P_0,P_4,P_6,\overline{P_7,P_8},P_9,P_{10}),\\
\displaystyle \ (P_0,P_4,\overline{P_6,P_8},P_9,P_{10}),
\ (P_0,\overline{P_4,P_8},P_9,P_{10}),\ (P_0,\overline{P_4,P_5},P_8,
P_9,P_{10})\ .
\end{array}$$

Consider the rational vector space
${\mathcal R}={\mathcal R}({\mathcal C})$ (only natural coordinates
will in fact occur) spanned by the set of all decorated 
roofs of a finite configuration ${\mathcal C}\subset {\mathbf R}^2$. 
Identifying a decorated
roof with the corresponding basis element, we define the
{\it transfer-matrix} $T$ as the linear application 
$T:{\mathcal R}\longrightarrow {\mathcal R}$ which sends
a decorated roof to the sum of its successors.
We introduce moreover the element $C_k\in{\mathcal R}(E)$
obtained by summing up all initially decorated roofs
of length $k$ with skyline $\partial_-{\mathcal C}$. 
Set $V_0=0$ and
$$V_i=C_i+T(V_{i-1})\ .$$
We denote by $W:{\mathcal R}\longrightarrow
{\mathbf Q}[s^{1/2}]$ the linear map which sends a decorated
roof to $0$ if its skyline is different from $\partial_+{\mathcal C}$
and wich sends a decorated roof $R$ with skyline $\partial_+{\mathcal
C}$ to $s^{\hbox{length}(R)/2}$.

\begin{prop} \label{ctpc}
The complete triangulation polynomial of a finite configuration
${\mathcal C}=(P_0,\dots,P_n)$ is given by
$$s\ \sum_{k=1}^n s^{k/2}\ W(V_k)\ .$$
\end{prop}

The proof is a straightforward consequence of Proposition \ref{maxroofseq}
and of the behaviour of lenghts of successive roofs.

An obvious modification of
the above proposition counts triangulations of a polygonal
(not necessarily convex) region squeezed between two skylines. 

If one is only interested in the number of maximal triangulations
of $\hbox{Conv}({\mathcal C})$, one considers the
transfer matrix $\tilde T$ which sends a decorated roof $R$ to 
the sum of its immediate successors. The vector $\tilde C_i$ is
zero except for $i+1=\sharp(\partial_-{\mathcal C}\cap {\mathcal C})$
and $\tilde C_i$ consists then of the initially decorated roof 
defined by the set $\partial_-{\mathcal C}\cap {\mathcal C}$.
We define $\tilde V_0=0$ and $\tilde V_i=\tilde C_i+\tilde T(\tilde V_{i-1})$.
Set now $\tilde W(R)=0$ except if $R$ is a decorated roof 
with underlying set $\partial_+{\mathcal C}\cap {\mathcal C}$
where we set $\tilde W(R)=1$. The number of maximal
triangulations of $\hbox{Conv}({\mathcal C})$ is now given by
$$\sum_i \tilde W(\tilde V_i)$$
(in fact, only one summand is non-zero).

Given a near-edge $E=(P_0,\dots,P_n)$, we introduce the vectors $
C_i\in{\mathcal R}(E)$ and the transfer matrix
$T\in\hbox{End}({\mathcal R}(E))$ as
above (where we consider $E$ as a configuration) and compute
the vectors $V_0=0$,
$$V_i=C_i+T(V_{i-1}),\ 1\leq i\leq n-1\ .$$

We introduce moreover the linear map $\tilde W:{\mathcal R}(E)\longrightarrow
{\mathbf Q}[s^{1/2},t]$ defined by
$$\begin{array}{lcl}
\displaystyle \tilde W(R)&=&s^{\hbox{length}(R)/2}\ p_{\hbox{length}(R)}\\
\displaystyle &=&s^{\hbox{length}(R)/2} \ \sum_k
(-1)^k{\hbox{length}(R)-k\choose k}\ t^{\hbox{length}(R)-k}
\end{array}$$
for any decorated roof $R$ of ${\mathcal C}$. One shows now easily
the following result.

\begin{prop} \label{cep} The complete edge-polynomial of a finite near-edge
$E=(P_0,\dots,P_n)$ is given by
$$\sum_{k=1}^n s^{k/2}\ \tilde W(V_k)\ .$$
\end{prop}

%%%%%%%%%%%%%%%%%%%%%%%%%%%%%%%%%%%%%%%%%%%%%%%%%%%%%%%%%%%%%%%%%%%%%%%%%%%

\subsection{A few algorithmic remarks}

Propositions \ref{ctpc}  and \ref{cep} yield easily algorithms for computing
complete triangulation and edge polynomials.

These algorithms are
of exponential complexity: properly organized their
memory and time requirements are of order $2^nf(n)$
where $f(n)$ is polynomial (of degree $1$ for memory requirements,
neglecting multiprecision problems,
and of fairly small degree for time requirements) for a
configuration of $n+1$ points. A nice feature of these algorithms is their
simplicity: As an example, we run them by hand 
for the number of triangulations squeezed between two
skylines of a configuration having $8$ points and  
for the near-edge $E_c$ of the near-polygon depicted in Figure 2. 

For computational purposes, it is useful to have a simple representation
for a decorated roof 
$$R=(P_{i_0}=P_0,P_{i_1},\dots,\overline{P_{i_d},P_{i_{d+1}}},\dots,
P_{i_{w-1}},P_{i_w}=P_n)$$
of a configuration (or near-edge) ${\mathcal C}=(P_0,\dots,P_n)$.
Such a roof $R$ is for instance encoded by the integer 
$l(R)\in\{0,\dots,n\cdot 2^{n-1}-1\}$
defined by
$$l(R)=d\cdot 2^{n-1}+\sum_{k=1}^{w-1}2^{i_k-1}\ .$$
This yields an injection from the set of decorated
roofs (of a given finite configuration $(P_0,\dots,P_n)$) 
into a subset
of $\sum_{k=1}^n k{n-1\choose k-1}=(n+1)2^{n-2}$ natural integers
$<n\cdot 2^{n-1}$ for $n>1$.
A natural integer $$m=\alpha\cdot 2^{n-1}+\sum_{i=0}^{n-2} \epsilon_i\ 2^i,\ 
0\leq \alpha\leq n-1,\ \epsilon_i\in\{0,1\}$$
is of the form $l(R)$ if and only if either $\alpha=0$ or
$\epsilon_{\alpha-1}=1$. The corresponding
roof $R$ is then given by
$$R=(P_0,\{P_{i+1}\ \vert\  0\leq i\leq n-2\hbox{ and }\epsilon_i=1\},P_n)$$
with decoration beginning at $P_{\alpha}$.

Let us also remark that many decorated roofs lead to dead ends
when computing triangulation polynomials. Indeed, consider a decorated
roof of the form
$$R=(P_{i_0},\dots,P_{i_k},\dots,\overline{P_{i_d},P_{i_{d+1}}},\dots)$$
where $P_{i_k}\in\partial_+{\mathcal C}$ for $k\leq d$. If the
partial skyline defined by the piece-wise linear path with
segments consecutive points of $(P_{i_0},\dots,P_{i_k})$
is not contained in $\partial_+{\mathcal C}$ then there exists
no roof-sequence containing $R$ which stops with a roof having skyline
$\partial_+({\mathcal C})$. 
Indeed, $R$ hits the \lq\lq ceiling'' $\partial_+({\mathcal C})$ at
$P_{i_k}$ before reaching the decorated edge. This
freezes the partial skyline $(P_{i_0},\dots,P_{i_k})$ in all successors
(and iterated successors) of $R$.

If one is only interested in computing the number of maximal triangulations
or the maximal edge-polynomial of a near-edge, a few simplifications take 
place: 

For non-generic configurations, only decorated roofs
which are {\it full} have to be considered where a roof is full if its
sequence contains all points of ${\mathcal C}$ lying on its skyline.

As already mentionned, the suitable
transfer matrix in this case is defined by considering only
immediate successors of a decorated roof (recall that a 
successor is immediate if it is obtained by adjoining a suitable
triangle which is minimal, i.e. which intersects ${\mathcal C}$ 
only in its three vertices).

Every decorated roof appears at most in one vector of
$V_1,\dots,V_n$. This makes memory requirements
somewhat smaller.
\bigskip

Let us also add that the algorithm which computes the complete
triangulation polynomial of $\hbox{Conv}({\mathcal C})$ 
can be modified for computing the
number of triangulations (with vertices in ${\mathcal C}$)
for any polygonal region enclosed by two suitable non-crossing
skylines $S_\alpha,\ S_\omega$: Replace $\partial_\pm({\mathcal C})$
by $S_\alpha$ and $S_\omega$.

The resulting 
algorithm for computing triangulation polynomials or edge-polynomials
is straightforward and best illustrated by examples. This
will be done in the 2 following subsections.

A few further tricks may be used:

Computing $V_i=C_i+T(V_{i-1})$ needs only $C_i$ and $V_{i-1}$.
This can be used to save memory by computing the contribution
$s^{(i+1)/2}W(V_{i-1})$ or $s^{(i-1)/2}\tilde W(V_{i-1})$ (in the case
of an edge polynomial)
of $V_{i-1}$ to the final result and by erasing $V_{i-1}$ after
determination of $V_i$.

Multiprecision problems can be avoided using the Chinese Remainder
Theorem.
Compute all cofficients of decorated roofs
modulo several primes and use them to reconstruct the result
(this needs however the computation of an approximation
of the final result which can for instance be done
by a computation using the type of floating reals).   

Decorated roofs should always be totally ordered
(using for instance the complete order coming from the 
labelling $l(R)$ for a decorated roof $R$). Since
there are at most $\sim n2^{n-2}$ such roofs, the coefficient of a
decorated roof can then be accessed using roughly
$O(n\ \hbox{log}(n))$ operations
using a classical divide and conquer stategy.

\subsection{A triangulation polynomial}

\begin{figure}[h]\label{expletetrys}
\epsfysize=2.5cm
\epsfbox{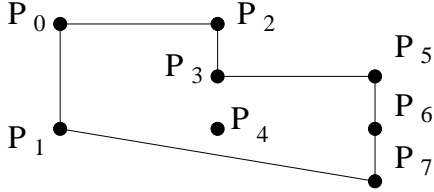}
\caption{A configuration enclosed by a polygonal line}
\end{figure}  
We want to compute the complete triangulation polynomial
of the (non-convex) polygonal region enclosed
by the two skylines $S_\alpha,S_\omega$ defined by
$S_\alpha=(P_0,P_1,P_7)$ and $S_\omega=(P_0,P_2,P_3,P_5,P_6,P_7)$
of the configuration
$$(P_0,\dots,P_7)=\left(\begin{array}{rrrrrrrrrrrrr}
0&0&1&1&1&2&2&2\\
3&1&3&2&1&2&1&0\end{array}\right)$$
depicted in Figure 13.

The following lists contain all relevant data. 
First we establish the list of all needed decorated roofs, enumerated 
and denoted using conventions and notations introduced in
the previous subsection. For concision,
we indicate a decorated roof $R=(P_0,\dots,\overline{P_{i_d},
P_{i_{d+1}}},\dots,P_7)$ by the (decorated) sequence
$(0\ i_1\ \cdots \overline{i_d\ i_{d+1}}\cdots 7)$
of the corresponding indices.
For every roof $R$, we indicate also its image $T(R)$
obtained by applying the transfer operator $T$,
except that roofs in $T(R)$ leading to dead-ends
(like $(01\overline{36}7)$ for instance) have been omitted.
For a roof $R$ with skyline $S_\omega$, we indicate $s^{\hbox{length}(R)/2}$
instead of $T(R)$ which is always zero in this case.

Data appearing in boldface can be omitted if one is only interested 
in the number of maximal triangulations.

$$\begin{array}{lll}
R_0&(\overline{07})&R_4\\
R_1&(\overline{01}7)&\mathbf{R_0+R_{69}+}R_{73}\mathbf{+R_{81}+R_{97}}\\
R_4&(\overline{03}7)&R_6 \\
R_6&(\overline{02}37)&\mathbf{ R_{150}}+R_{166} \\

R_8&(\overline{04}7)&R_0+R_{12} \\
R_9&(\overline{01}47)&R_8+R_{77}\mathbf{+R_{153}}+R_{169} \\
R_{12}&(\overline{03}47)&R_{14} \\
R_{14}&(\overline{02}347)&R_{134}+\mathbf{ R_{222}}+R_{238} \\
\mathbf{ R_{20}}&\mathbf{ (\overline{03}57)}&\mathbf{ R_{22}} \\
\mathbf{ R_{22}}&\mathbf{ (\overline{02}357)}&\mathbf{ s^2} \\
R_{36}&(\overline{03}67)&R_{38} \\
R_{38}&(\overline{02}367)&R_{182} \\
R_{52}&(\overline{03}567)&R_{54} \\
R_{54}&(\overline{02}3567)&s^{5/2} \\
\end{array}$$

$$\begin{array}{lll}
R_{69}&(0\overline{13}7)&R_4+R_{71} \\
R_{71}&(0\overline{12}37)&R_6 \\
R_{73}&(0\overline{14}7)&R_8+R_{77}+\mathbf{R_{153}}+R_{169} \\
R_{77}&(0\overline{13}47)&R_{12}+R_{79} \\
R_{79}&(0\overline{12}347)&R_{14}\\
\mathbf{ R_{81}}&\mathbf{ (0\overline{15}7)}&\mathbf{ R_{85}} \\
\mathbf{ R_{85}}&\mathbf{ (0\overline{13}57)}&\mathbf{ R_{20}+R_{87}}\\
\mathbf{ R_{87}} &\mathbf{ (0\overline{12}357)}&\mathbf{ R_{22}} \\
\mathbf{R_{97}}&\mathbf{(0\overline{16}7)}&\mathbf{R_{101}+R_{113}} \\
R_{101}&(0\overline{13}67)&R_{36}+R_{103}\\
R_{103}&(0\overline{12}367)&R_{38} \\
R_{113}&(0\overline{15}67)& R_{117}\\
R_{117}&(0\overline{13}567)&R_{52}+R_{119}\\
R_{119}&(0\overline{12}3567)&R_{54}
\end{array}$$

$$\begin{array}{lll}
R_{134}&(02\overline{37})&\mathbf{ R_{150}}+R_{166} \\
\mathbf{ R_{150}}&\mathbf{ (02\overline{35}7)}&\mathbf{  s^2}\\
\mathbf{ R_{153}}&\mathbf{ (01\overline{45}7)}&\mathbf{  R_{81}}\\
R_{166}&(02\overline{36}7)& R_{182}\\
R_{169}&(01\overline{46}7)& R_{185}\\
R_{182}&(02\overline{35}67)&s^{5/2}\\
R_{185}&(01\overline{45}67)& R_{113}\\
\end{array}$$

$$\begin{array}{lll}
\mathbf{ R_{222}}&\mathbf{ (023\overline{45}7)}&\mathbf{ R_{150}}\\
R_{238}&(023\overline{46}7)&R_{166}+R_{254}\\
R_{254}&(023\overline{45}67)&R_{182}
\end{array}$$

The vectors $C_i$ are all $0$ except $C_2=R_1$. We have then $V_1=0$,
$V_2=C_2=R_1$ and $V_i=T(V_{i-1}),\ i=3,\dots,9$ where $T$ is the transfer
operator.

$$\begin{array}{lcl}
V_2&=&C_2=R_1\\
V_3&=&\mathbf{R_0+R_{69}+}R_{73}\mathbf{+R_{81}+R_{97}}\\
V_4&=&\mathbf{R_4+(R_4+R_{71})+}(R_8+R_{77}\mathbf{+R_{153}}+R_{169})\mathbf{+R_{85}+(R_{101}+R_{113})}\\
&=&\mathbf{2R_4+}R_8\mathbf{+R_{71}+}R_{77}\mathbf{+R_{85}+R_{101}+R_{113}+R_{153}}+R_{169}\\
V_5&=&\mathbf{2R_6+}(R_0+R_{12})\mathbf{+R_6}+(R_{12}+R_{79})\mathbf{
+(R_{20}+R_{87})}\\
&&\quad \mathbf{+(R_{36}+R_{103})+R_{117}+R_{81}}+R_{185}\\
&=&R_0\mathbf{+3R_6+}2R_{12}\mathbf{+R_{20}+R_{36}}+R_{79}\mathbf{+R_{81}}\\
&&\quad \mathbf{+R_{87}+R_{103}+R_{117}}+R_{185}\\
V_6&=&R_4\mathbf{+3(R_{150}+R_{166})}+2R_{14}\mathbf{+R_{22}+R_{38}}+R_{14}
\mathbf{+R_{85}}\\
&&\quad \mathbf{+R_{22}+R_{38}+(R_{52}+R_{119})}+R_{113}\\
&=&R_4+3R_{14}\mathbf{+2R_{22}+2R_{38}}\\
&&\quad \mathbf{+R_{52}+R_{85}}+R_{113}\mathbf{+R_{119}+3R_{150}+3R_{166}}\\
V_7&=&R_6+3(R_{134}\mathbf{+R_{222}}+R_{238})\mathbf{+2\cdot 0+2R_{182}}\\
&&\quad \mathbf{+R_{54}+(R_{20}+R_{87})}+R_{117}\mathbf{+R_{54}+3\cdot 0
+3R_{182}}\\
&=&R_6\mathbf{+R_{20}+2R_{54}+R_{87}}+R_{117}\\
&&\quad +3R_{134}\mathbf{+5R_{182}+3R_{222}}+3R_{238}\\
V_8&=&(\mathbf{R_{150}+}R_{166})\mathbf{+R_{22}+2\cdot 0+R_{22}}+(R_{52}+R_{119})\\
&&\quad +3(\mathbf{R_{150}+}R_{166})\mathbf{+5\cdot 0+3R_{150}}+
3(R_{166}+R_{254})\\
&=&\mathbf{2R_{22}+}R_{52}+R_{119}\mathbf{+7R_{150}}+7R_{166}+3R_{254}\\
V_9&=&\mathbf{2\cdot 0+}R_{54}+R_{54}\mathbf{+7\cdot 0}+7R_{182}+3R_{182}\\
&=&2R_{54}+10R_{182}\\
V_{10}&=&0
\end{array}$$

We have $W(V_2)=\dots=W(V_5)=0$ and
$$
(W(V_6),\dots,W(V_9))=(5s^2,7s^{5/2},9s^2,12s^{7/2})$$
and the complete triangulation polynomial of the polygonal region 
between $S_\alpha$ and $S_\omega$ is given by
$$s(s^3\ 5s^2+s^{7/2}\ 7s^{5/2}+s^4\ 9s^2+s^{9/2}\ 12s^{7/2})=
12s^8+16s^7+5s^6\ .$$

\subsection{An edge-polynomial}

\begin{figure}[h]\label{edgeE5c}
\epsfysize=3cm
\epsfbox{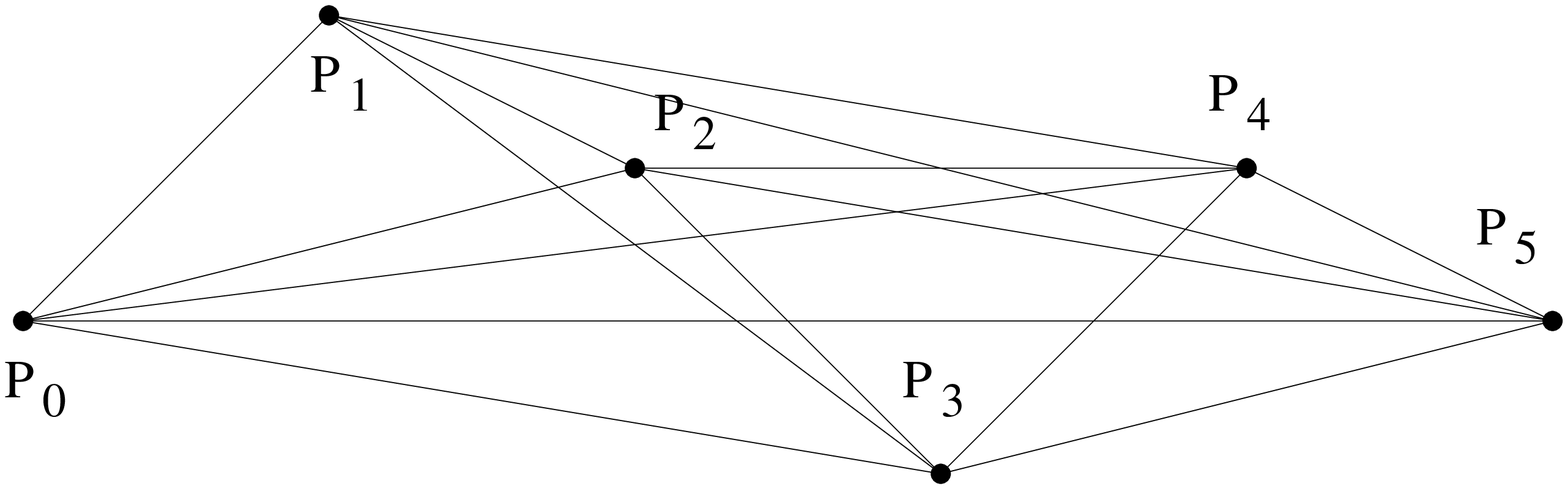}
\caption{The near-edge $E_c$}
\end{figure}  

We consider the $5-$near-edge $E=E_c$ represented in Figure 14
of the near-polygon appearing in Figure 2. This near-edge $E=(P_0,\dots, P_5)$
can be represented by
$$E=\left(\begin{array}{rrrrrr}0&1&2&3&4&5\\0&2&1&-1&1&0\end{array}\right)$$
In the following table we list all decorated roofs
arising during the computations. 
The first two columns list all $23$ decorated roofs needed for 
the computation, the third column
yields the image $T(R_i)$ (where $T$ denotes the transfer matrix)
of the decorated roof
$R_i$ described by the previous columns. The $6-$th line for instance
lists the roof $(\overline{P_0,P_1},P_3,P_5)$ encoded
by $(\overline{01}35)$ and wich we call $R_5$ using the labelling
introduced previously.
An easy computation shows that 
$$T(R_5)=(0\overline{15})+(0\overline{12}35)+(01\overline{34}5)\ .$$
Using our numerotation of decorated roofs, we have thus
$$T(R_5)=R_{17}+R_{23}+R_{45}\ .$$
Data appearing in boldface is only involved in the computation
of the complete edge-polynomial and can be omitted for the 
computation of the maximal edge-polynomial.

$$\begin{array}{llll}
R_0&(\overline{05})& {\bf (\overline{01}5)+}(\overline{02}5)+
(\overline{04}5)&  {\bf R_1+}R_2+R_8\\
{\bf R_1}&{\bf (\overline{01}5)}&  (0\overline{14}5)&R_{25} \\
R_2&(\overline{02}5)&  (\overline{01}25)+(0\overline{24}5)&R_3+R_{26} \\
R_3&(\overline{01}25) & (0\overline{15})+(01\overline{24}5)&R_{17}+R_{43}\\
R_4&(\overline {03}5)& (\overline{05})+(\overline{01}35)+(\overline{02}35)
+(0\overline{34}5)&R_0+R_5+R_6+R_{28}\\
R_5&(\overline{01}35)& {\bf (0\overline{15})+}(0\overline{12}35)+
(01\overline{34}5)&{\bf R_{17}+}R_{23}+R_{45}\\
R_6&(\overline{02}35)& (0\overline{25})+(\overline{01}235)+(02\overline{34}5)&R_{18}+R_7+R_{46}\\
R_7&(\overline{01}235) & (01\overline{25})+(012\overline{34}5)&R_{35}+R_{63}\\
R_8&(\overline{04}5)& {\bf (\overline{01}45)+}(\overline{02}45)&{\bf R_9+}R_{10}  \\
{\bf R_{9}}&{\bf (\overline{01}45)} & 0&0\\
R_{10}&(\overline{02}45) & (\overline{01}245)& R_{11}\\
R_{11}&(\overline{01}245)&(0\overline{14}5)& R_{25}\\
R_{17}&(0\overline{15})& (0\overline{14}5)&R_{25}  \\
R_{18}&(0\overline{25}) & (0\overline{24}5)&R_{26}\\
R_{23}&(0\overline{12}35) & (01\overline{25})+(012\overline{34}5)&R_{35}+
R_{63}\\
R_{25}&(0\overline{14}5) & 0&0\\
R_{26}&(0\overline{24}5) & 0&0\\
R_{28}&(0\overline{34}5)& (\overline{04}5)&R_8 \\
R_{35}&(01\overline{25}) & (0\overline{15})+(01\overline{24}5)&
R_{17}+R_{43}\\
R_{43}&(01\overline{24}5)& (0\overline{14}5)& R_{25}\\
R_{45}&(01\overline{34}5) & {\bf (0\overline{14}5)}&{\bf R_{25}}\\
R_{46}&(02\overline{34}5) & (0\overline{24}5)&R_{26}\\
R_{63}&(012\overline{34}5)& (01\overline{24}5) & R_{43}\\
\end{array}
$$

The vectors $C_i$ are $0$ except $C_2=R_4$.

\begin{eqnarray*}
V_1&=&C_1=0\\
V_2&=&T(V_1)+C_2=C_2=R_4\\
V_3&=&T(V_2)+C_3=R_0+R_5+R_6+R_{28}\\
V_4&=&T(V_3)+C_4=({\bf R_1+}R_2+R_8)+({\bf R_{17}+}R_{23}+R_{45})\\
&&\quad +(R_{18}+R_7+R_{46})+R_8\\
&=&{\bf R_1+}R_2+R_7+2R_8{\bf +R_{17}}+R_{18}+R_{23}+R_{45}+R_{46}\\
V_5&=&T(V_4)+C_5={\bf R_{25}+}(R_3+R_{26})+(R_{35}+R_{63})
+2({\bf R_9+}R_{10})\\&&\quad {\bf +R_{25}}
+R_{26}+(R_{35}+R_{63}){\bf +R_{25}}+R_{26}\\
&=&R_3{\bf +2R_9}+2R_{10}{\bf +3R_{25}}+3R_{26}+2R_{35}+2R_{63}\\
V_6&=&T(V_5)+C_6=(R_{17}+R_{43}){\bf +2\cdot 0}+2R_{11}{\bf +3\cdot 0}
+3\cdot 0\\
&&\quad+2(R_{17}+R_{43})+2R_{43}\\
&=&2R_{11}+3R_{17}+5R_{43}\\
V_7&=&T(V_6)+C_7=3R_{25}+5R_{25}+2R_{25}\\
&=&10R_{25}
\end{eqnarray*}

\begin{eqnarray*}
{\overline{p}(E)}&=&\sum_i s^{i/2}\tilde W(V_i)\\
&=&p_2s^2+(p_1s^2+3p_3s^3)+(6p_2s^3+4p_4s^4)\\
&&\quad +(13p_3s^4+2p_5s^5)+(3p_2s^4+7p_4s^5)+10p_3s^5\\
&=&(10p_3+7p_4+2p_5)s^5+(3p_2+13p_3+4p_4)s^4\\
&&\quad +(6p_2+3p_3)s^3+(p_1+p_2)s^2
\end{eqnarray*}

%%%%%%%%%%%%%%%%%%%%%%%%%%%%%%%%%%%%%%%%%%%%%%%%%%%%%%%%%%%%%%%%%%%%%%%

Roland Bacher, Institut Fourier, UMR 5582,
Laboratoire de Math\'ematiques, BP 74, 38402 St. Martin d'H\`eres Cedex,
France, Roland.Bacher@ujf-grenoble.fr

\end{document}